\documentclass[10pt,letterpaper,twocolumn]{article}
\usepackage[left=0.65in,right=0.65in,top=0.8in,bottom=0.8in]{geometry}
\usepackage{cite}
\usepackage{amsmath,amssymb,amsfonts}
\usepackage{amsthm}
\usepackage{graphicx}
\usepackage{algorithm,algorithmic}
\usepackage{hyperref}
\hypersetup{hidelinks}
\usepackage{textcomp}
\usepackage{mathrsfs}
\usepackage{enumerate}
\usepackage{float}
\usepackage{tikz}
\usepackage{etoolbox}
\makeatletter

\usepackage{titlesec}
\titleformat{name=\section,numberless}[block]
  {\normalfont\normalsize\bfseries\centering}
  {}{0pt}{}
\titleformat{\section}
  {\normalfont\normalsize\bfseries\centering}
  {\thesection}{0.6em}{}

\titleformat{\subsection}
  {\normalfont\normalsize\itshape\raggedright}
  {\thesubsection}{0.6em}{}

\titleformat{\subsubsection}
  {\normalfont\small\itshape\raggedright}
  {\thesubsubsection}{0.6em}{}

\titlespacing*{\section}{0pt}{1.2ex plus .3ex}{0.6ex}
\titlespacing*{\subsection}{0pt}{1.0ex plus .3ex}{0.4ex}
\titlespacing*{\subsubsection}{0pt}{0.8ex plus .2ex}{0.3ex}

\usetikzlibrary{
    arrows.meta,
    positioning,
    fit,
    decorations.markings,
    calc,
    intersections
}

\allowdisplaybreaks[4]

\tikzset{
    degil/.style={
        postaction={
            decorate,
            decoration={
                markings,
                mark=at position 0.5 with {
                    \draw[-] (-2pt,-2pt) -- (2pt,2pt);
                }
            }
        }
    }
}

\newtheorem{theorem}{Theorem}[section]
\newtheorem{lemma}[theorem]{Lemma}
\newtheorem{corollary}[theorem]{Corollary}
\newtheorem{proposition}[theorem]{Proposition}

\newtheorem{definition}[theorem]{Definition}

\newtheorem{remark}[theorem]{Remark}

\newtheorem{assumption}{Assumption}

\AtEndEnvironment{example}{\hfill\ensuremath{\blacksquare}}

\makeatletter
\DeclareOldFontCommand{\rm}{\normalfont\rmfamily}{\mathrm}
\DeclareOldFontCommand{\sf}{\normalfont\sffamily}{\mathsf}
\DeclareOldFontCommand{\tt}{\normalfont\ttfamily}{\mathtt}
\DeclareOldFontCommand{\bf}{\normalfont\bfseries}{\mathbf}
\DeclareOldFontCommand{\it}{\normalfont\itshape}{\mathit}
\DeclareOldFontCommand{\sl}{\normalfont\slshape}{\@nomath\sl}
\DeclareOldFontCommand{\sc}{\normalfont\scshape}{\@nomath\sc}
\makeatother

\usepackage{csquotes}  
\newcommand\q{\enquote}

\usepackage{physics}   

\newcommand \N   {\mathbb{N}}
\newcommand \R   {\mathbb{R}}

\newcommand \C   {\mathbb{C}}

\newcommand \Kinf{\mathcal{K_\infty}}

\newcommand \KL  {\mathcal{KL}}

\newcommand{\Dc}{\ensuremath{\mathcal{D}}}

\newcommand{\vertiii}[1]{{\left\vert\kern-0.25ex\left\vert\kern-0.25ex\left\vert #1 
    \right\vert\kern-0.25ex\right\vert\kern-0.25ex\right\vert}}

\newcommand \qrq   {\quad\Rightarrow\quad}

\newcommand \sis   {\ \ \Iff\ \ }
\newcommand \srs   {\ \ \Rightarrow\ \ }

\newcommand \Iff   {\Leftrightarrow}

\newcommand{\normt}[1]{{\left\vert\kern-0.25ex\left\vert\kern-0.25ex\left\vert #1 
		\right\vert\kern-0.25ex\right\vert\kern-0.25ex\right\vert}}

\newif\ifMath					
\newif\ifEngi					

\newif\ifDFGtext					 

\newif\ifAndo              
													
\newif\ifExercises					
\newif\ifSolutions          
\newif\ifGerman							
\newif\ifEnglish						

\newif\ifnothabil						

\newif\ifFuture							

\newif\ifConf                    
\newif\ifJournal								 

\newif\ifNOTFORBOOK
\newif\ifFullVersion
\newif\ifExludedDueToSpaceReasons

\usepackage{xifthen}

\newcommand{\einsnorm}[2]{\ensuremath{
    \!\!\;\!\!\!\;
    \left\bracevert\!\!\!\!\!\left\bracevert
    \!
		\ifthenelse{\isempty{#2}}{#1}{#1(#2)}
    \!
      \right\bracevert\!\!\!\!\!\right\bracevert
    \!\!\;\!\!\!\;
  }}

\usepackage{xcolor}
\definecolor{blond}{rgb}{0.98, 0.94, 0.75}
	
\newlength\mytemplen
\newsavebox\mytempbox

\makeatletter
\newcommand\mybluebox{%
    \@ifnextchar[
       {\@mybluebox}%
       {\@mybluebox[0pt]}}

\def\@mybluebox[#1]{%
    \@ifnextchar[
       {\@@mybluebox[#1]}%
       {\@@mybluebox[#1][0pt]}}

\def\@@mybluebox[#1][#2]#3{
    \sbox\mytempbox{#3}%
    \mytemplen\ht\mytempbox
    \advance\mytemplen #1\relax
    \ht\mytempbox\mytemplen
    \mytemplen\dp\mytempbox
    \advance\mytemplen #2\relax
    \dp\mytempbox\mytemplen
    \colorbox{blond}{\hspace{1em}\usebox{\mytempbox}\hspace{1em}}}

\makeatother

\makeatletter
\let\origd=\d
\renewcommand*\d{
  \relax\ifmmode
    \mathrm{d}%
  \else
    \expandafter\origd
  \fi
}\makeatother

\usepackage{mathtools}

\makeatletter 
\newcommand{\pushright}[1]{\ifmeasuring@#1\else\omit\hfill$\displaystyle#1$\fi\ignorespaces}
\newcommand{\pushleft}[1]{\ifmeasuring@#1\else\omit$\displaystyle#1$\hfill\fi\ignorespaces}
\makeatother 

\newcounter{syscounter}

\newcounter{WPcounter}
\newcounter{PRcounter}
\def\BibTeX{{\rm B\kern-.05em{\sc i\kern-.025em b}\kern-.08em
    T\kern-.1667em\lower.7ex\hbox{E}\kern-.125emX}}

\begin{document}
\title{Lyapunov Theory for Output-to-State Stability of Nonlinear Infinite-Dimensional Systems%
\thanks{This work was partially supported by the German Research Foundation
(DFG) under Grants MI 1886/2-2 (Project No. 415101813) and MI 1886/3-1
(Project No. 540580186).
\protect\endgraf
Q. Chen and F. Wirth are with the Faculty of Computer Science and
Mathematics, University of Passau, 94032 Passau, Germany
(e-mail: \protect\nolinkurl{qiaoling.chen@uni-passau.de};
\protect\nolinkurl{fabian.lastname@uni-passau.de}).
\protect\endgraf
A. Mironchenko is with the Department of Mathematics, University of
Bayreuth, 95447 Bayreuth, Germany
(e-mail: \protect\nolinkurl{andrii.mironchenko@uni-bayreuth.de}).}}

\author{%
  Qiaoling Chen,
  Andrii Mironchenko, and
  Fabian Wirth
}

\date{}
\maketitle

\date{}\maketitle

\begin{abstract}
For a broad class of nonlinear infinite-dimensional systems with outputs, we show that
output-to-state stability (OSS) is equivalent to 
the existence of a coercive OSS Lyapunov function in implication form. 
While establishing this equivalence, we prove that OSS is equivalent to uniform global asymptotic stability modulo output, relate these properties to the vanishing output vanishing state property, and highlight new challenges that do not arise in the nonlinear finite-dimensional or linear infinite-dimensional settings. The results are illustrated by a nonlinear reaction-diffusion system with a nonlinear output map.
\end{abstract}

\par

\section{Introduction}

The notion of detectability, introduced in \cite{wonham1968matrix} for linear systems, is paramount for observer design \cite[Section~8.1]{Sontag2008} and output-feedback stabilization \cite{logemann1992stability}. Detectability also plays an important role in nonlinear stability analysis \cite{Jiang1994} and state-feedback design, including LaSalle’s invariance principle and Jurdjevic–Quinn control \cite{angeli2004separation}. 
A linear system is called detectable if any two trajectories producing the same output converge to each other. A related notion is zero-detectability, which requires trajectories with identically zero output to converge to the origin. For linear ordinary differential equations (ODEs), detectability and zero-detectability are equivalent.

\textbf{Output-to-state stability.} 
Beyond convergence of zero-output trajectories, it is of interest to have estimates of the state magnitude in terms of the measured output that also capture transient behavior.
These state estimates are provided by  output-to-state stability (OSS)  introduced in \cite{sontag1996}.
 OSS is defined by estimating the norm of the trajectory in terms of the sum of a decaying transient and a term depending on the output. 
For linear ODEs without inputs, OSS is equivalent to detectability \cite{sontag1997output,krichman2001input}. This makes it a natural candidate for extending classical detectability to nonlinear systems, as proposed in \cite{sontag1996}. OSS estimates can be obtained using a function that decreases along trajectories whenever the state is sufficiently large compared to the output. Similar conditions on storage functions had already been used in earlier formulations of nonlinear detectability \cite[Remark~5]{sontag1997output}.  For nonlinear ODEs, OSS has been characterized in terms of OSS Lyapunov functions in dissipation form in \cite{sontag1996}.

The success of OSS theory led to the introduction of a variety of concepts that refine or extend it.
For ODE systems with inputs and outputs, input/output-to-state stability (IOSS) has been introduced in \cite{sontag1997output,krichman2001input} to formalize robust detectability of nonlinear systems. IOSS is equivalent to the existence of a smooth IOSS Lyapunov function \cite{krichman2001input}.
To obtain a proper extension not only of zero-detectability but also of detectability, the concept of incremental input/output-to-state stability was introduced in 
\cite{sontag1997output}.

The applications of (I)OSS theory are manifold. An important application is state estimation.  For nonlinear ODEs, global complete observability together with OSS guarantees the existence of a globally convergent observer \cite{astolfi2006global}. Incremental IOSS is a necessary condition for the existence of a robust full-state observer \cite{sontag1997output} and 
other types of state estimators, such as the extended Kalman filter, full information estimation, and moving horizon estimation; see \cite{allan2021nonlinear}. 
Incremental IOSS is widely used for  optimization-based state estimation \cite{RMD17}, and its relations to incremental dissipativity were studied in \cite{SHB12}. 
For optimal control problems with nonnegative stage cost on bounded sets, the existence of an IOSS-like Lyapunov function implies strict dissipativity for stage cost depending on the norms of the output and input \cite[p. 459 \& Corollary 1]{HoG19}.
In turn, strict dissipativity implies the turnpike property that ensures similarity between optimal trajectories on varying time horizons \cite{GrM16, MGA15, TrZ15}. This is essential for model predictive control schemes \cite{GrP17,RMD17,AlZ00}.

The fundamental idea of OSS is similar to that of
 input-to-state stability (ISS). The ISS framework has revolutionized the robust stability analysis of nonlinear control systems and their interconnections \cite{sontag1989smooth, Jiang1994, Sontag1996b, Dashkovskiy2010}, see the monograph \cite{mironchenko2023} for a comprehensive treatment of the ISS theory of ODEs.
ISS quantifies the input-state behavior of a system. It ensures that bounded inputs lead to bounded states and that the state converges to zero if the input decays to zero. ISS can be characterized as a combination of IOSS and input-to-output stability (IOS), see \cite[Proposition 3.1]{Jiang1994} for the ODE case and \cite[Section V.B]{Bachmann2024} for general infinite-dimensional systems.

These and many other contributions constitute comprehensive theories of output-to-state and input-to-state stability for ODE systems. 
In the last decade, the input-to-state stability framework for infinite-dimensional systems has been rapidly developed \cite{jayawardhana2008infinite,Karafyllis:2016a, Karafyllis:2018iss,Jacob:2018_SIAM,mironchenko2017characterizations,mironchenko2020input,Schwenninger2020,zheng2021approximations}, to cite a few. 
In contrast, (I)OSS theory remains largely unexplored in the infinite-dimensional setting.
 Sufficient conditions for IOSS and integral IOSS have been established for nonlinear impulsive delay systems, including systems with switching \cite{jiang2019input,hong2019input}.

As for OSS, in \cite{chen2026detectability}, we have initiated the development of the exponential output-to-state stability (eOSS) theory for linear infinite-dimensional systems with bounded output operators in Banach spaces. The main result of \cite{chen2026detectability} shows that eOSS is equivalent to the existence of a coercive eOSS Lyapunov function in implication form. 
Furthermore, it is shown that exponential detectability is strictly stronger than the existence of a coercive eOSS Lyapunov function in dissipation form and therefore eOSS, which is strictly stronger than zero-detectability. 
The linear theory developed in \cite{chen2026detectability} reveals a fundamental difficulty that persists in the nonlinear setting. The counterexample in \cite[Proposition 5.2]{chen2026detectability} shows that already for linear infinite-dimensional systems OSS does not, in general, imply the existence of an OSS Lyapunov function in dissipation form, in striking contrast to the ODE case \cite{sontag1996}.  This raises the question of whether OSS for nonlinear infinite-dimensional systems can be characterized by  OSS Lyapunov functions in implication form.

Transitioning from linear to nonlinear infinite-dimensional systems brings further difficulties. The construction of an exponential OSS Lyapunov function for linear OSS systems developed in \cite{chen2026detectability} relies on the linear semigroup structure and on an equivalent norm construction. These tools cannot be directly applied to the nonlinear dynamical systems considered here. At the same time, in contrast to the nonlinear ODE setting, infinite-dimensional systems may lack differentiability of the flow with respect to the initial state and smoothing techniques available in Euclidean spaces, which are crucially important for converse OSS Lyapunov constructions in \cite{sontag1996}.

Two related properties provide further perspectives on OSS. Uniform global asymptotic stability modulo output (UGASMO), introduced in \cite[Definition~2.9]{krichman2001input}, requires a uniform decay estimate for the state as long as its magnitude remains sufficiently large compared to the output.  The vanishing output vanishing state (VOVS) property requires the state to converge to zero whenever the output converges to zero.  It has recently been used as a sufficient condition for asymptotic convergence of nonlinear closed-loop systems \cite{preuster2026stabilization}. VOVS can also be viewed as a strengthening of zero-detectability by allowing outputs that merely converge to zero.  It is related in spirit to the converging-input converging-state (CICS) property in \cite{Sontag1989}. Necessary and sufficient conditions for the CICS property were later derived for a class of Lur’e systems in \cite{BillEtAl2017}.

\textbf{Contribution.} 
This paper develops a Lyapunov theory for OSS of a general class of nonlinear dynamical systems.

\begin{figure}[htbp]    
\vspace{.3cm}    
\centering
\begin{tikzpicture}
\node (BlockEquiv) at (-6,0) {OSS};

\node[above = 0.7cm of BlockEquiv] (e2) {$\exists$ coercive OSS Lyapunov function in implication form}; 


\node[above = 0.7cm of e2](e1)
{$\exists$ coercive OSS Lyapunov function in dissipation form};


 \node[below = 0.7cm of BlockEquiv] (h) {UGASMO};

 \node[below = 0.7cm of h] (g) {VOVS};

 \node[below = 0.7cm of g] (g1) {zero-detectability};

\path([xshift=-0.2cm]e2.south) edge[thick,double,double equal sign distance,-{Implies[]}] node[left=.2cm]{\footnotesize Theorem~\ref{thm:abdirect OSS LT}} ([xshift=-0.2cm] BlockEquiv.north)

([xshift=0.2cm]BlockEquiv.north) edge[thick,double,double equal sign distance,-{Implies[]}] node[right=.2cm]{\footnotesize Theorem~\ref{prop:E0V1} (under BIC)} ([xshift=0.2cm]e2.south);

\path([xshift=-0cm]BlockEquiv.south) edge[thick,double,double equal sign distance,{Implies[]}-{Implies[]}] node[right=.2cm]{\footnotesize Theorem~\ref{thm:oss-equi-ugasmo}} ([xshift=-0cm] h.north);

\path([xshift=-0.2cm]e1.south) edge[thick,double,double equal sign distance,-{Implies[]}] node[left=.2cm]{\footnotesize Proposition~\ref{prop:dflf-criterion}} ([xshift=-0.2cm] e2.north)

([xshift=0.2cm]e2.north) edge[thick,double,double equal sign distance,-{Implies[]},degil,draw=red] node[right=.2cm]{\footnotesize Remark~\ref{re:im-not-di}} ([xshift=0.2cm]e1.south)

([xshift=-0.2cm]h.south) edge[thick,double,double equal sign distance,-{Implies[]}] node[left=.2cm]{\footnotesize  Proposition~\ref{prop:UGASMO-implies-weak-convergence} (under BIC)} ([xshift=-0.2cm] g.north)

([xshift=0.2cm]g.north) edge[thick,double,double equal sign distance,-{Implies[]},degil,draw=red] node[right=.2cm]{\footnotesize Proposition \ref{prop:vovs-not-ugasmo}} ([xshift=0.2cm]h.south)

([xshift=0.2cm]g1.north) edge[thick,double,double equal sign distance,-{Implies[]},degil,draw=red] node[right=.2cm]{\footnotesize Proposition~\ref{prop:zero-dete-not-Vovs}} ([xshift=0.2cm]g.south)

([xshift=-0.2cm]g.south) edge[thick,double,double equal sign distance,-{Implies[]}] node[left=.2cm]{\footnotesize Remark~\ref{re:vovs-zero-dete}} ([xshift=-0.2cm]g1.north);

\end{tikzpicture}
\caption{Summary of the implications and non-implications.  } 
\label{fig:Result}
\end{figure}
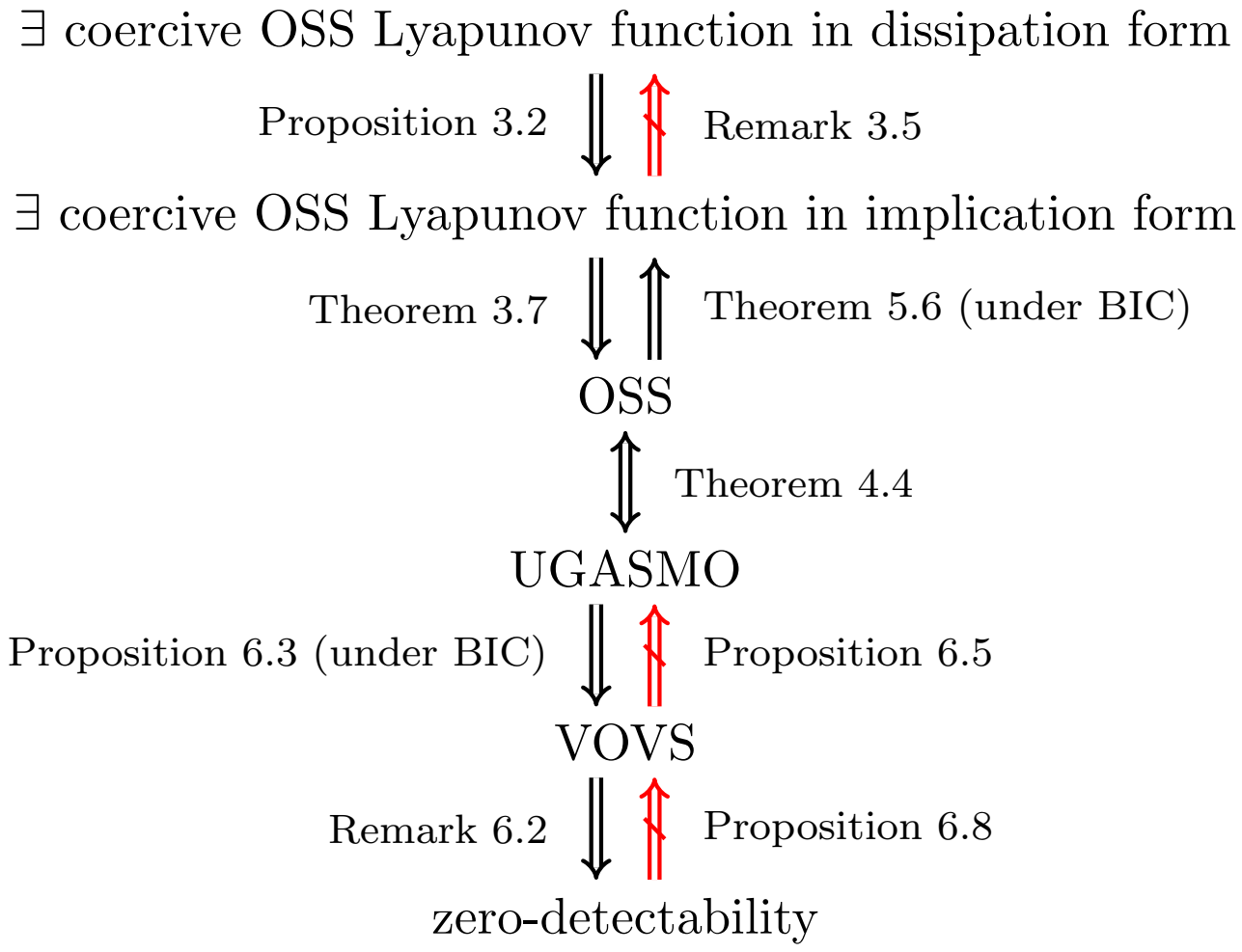

The main results and their logical relations are summarized in Figure~\ref{fig:Result}.
We first establish a direct OSS Lyapunov theorem.
Next, we establish the equivalence between OSS and UGASMO.  This equivalence is shown in the central part of Figure~\ref{fig:Result} and is then used to derive our main converse result:
\emph{Under the boundedness-implies-continuation (BIC) property, OSS, UGASMO, and the existence of a coercive OSS Lyapunov function in implication form are equivalent.}

As we have mentioned above and illustrated in the upper part of Figure~\ref{fig:Result}, already in the linear infinite-dimensional setting, the equivalence between the existence of OSS Lyapunov functions in implication and dissipation form breaks down \cite{chen2026detectability}.

To prove the converse Lyapunov theorem, we use the equivalence between OSS and UGASMO and construct a Lyapunov-like function directly along trajectories as a weighted supremum taken before the first entrance into an output-dominated region. The semigroup property and the shift property of the entrance time yield the required decay estimate, and a continuous cutoff is used to extend the function to the whole state space.

As summarized in the lower part of Figure~\ref{fig:Result}, under the BIC property, UGASMO is shown to imply VOVS, while VOVS implies zero-detectability. Two counterexamples demonstrate that the corresponding converse implications fail in general.

Finally, the developed theory is applied to a nonlinear reaction-diffusion system with a nonlinear output map, for which an explicit coercive OSS Lyapunov function is constructed.

A conference version of this paper, which does not contain the converse Lyapunov results, was presented in \cite{chen2025lyapunov}.

\textbf{Outlook.} This paper can serve as a starting point for a broad research program for the development of the foundations of the OSS theory. The next step could be the development of (incremental) IOSS theory, relating it to the strict dissipativity, turnpike properties, and observer design of nonlinear distributed parameter systems.

\textbf{Organization.} The remainder of the paper is organized as follows. Section~\ref{sec:pro setting} introduces the dynamical system framework and the basic notions used throughout the paper. Section~\ref{sec:direct OSS} develops direct OSS Lyapunov results and studies the relation between implication and dissipation forms. Section~\ref{sec:gasmo} introduces UGASMO and establishes its equivalence with OSS. Section~\ref{sec:converse OSS} is devoted to the converse OSS Lyapunov theorem. Section~\ref{sec:vovs} investigates the VOVS property and its relations to UGASMO and zero-detectability. Section~\ref{sec:example} applies the developed theory to a nonlinear reaction-diffusion system with a nonlinear output map. Section~\ref{sec:discuss} discusses the results and their relation to existing work.  Auxiliary technical results are collected in the Appendix.

\textbf{Notation.}  Let $\mathbb{R}$ be the set of real numbers, $\mathbb{R}_+:=[0,\infty)$ be the set of nonnegative real numbers and $\mathbb{N}$ be the set of positive integers. We write $\C$ for the complex numbers and $i$ for the imaginary unit.  Given topological spaces $X, Y$, the  space of continuous functions from $X$ to $Y$ is denoted by
$C(X, Y)$. We denote for short $C(X):=C(X,\R)$. For an interval $I\subset \R$ and $k\in \N$, $C^k(I)$
 denotes the space of real-valued functions that are $k$-times continuously differentiable. We use $C^k(I,Z)$ for the corresponding space of $Z$-valued functions, where $Z$ is a normed linear space; derivatives at included endpoints are one-sided. For normed spaces $X, Y,$ $\mathcal{L}(X,Y)$ denotes the bounded linear operators from $X$ to $Y$, and $\mathcal{L}(X):=\mathcal{L}(X,X).$

We use the following classes of  comparison functions \cite{Kel14}:
\begin{align*}
\mathcal{P}&:=\{\gamma:\mathbb{R}_+\rightarrow\mathbb{R}_+\mid\gamma\ \text{is continuous}\ \text{with}\ \gamma(0)=0\\
&\qquad\qquad\qquad \qquad \quad\text{and}\ \gamma(r)>0\ \text{for}\ r>0\},\\
\mathcal{K}&:=\{\gamma\in\mathcal{P}\mid \gamma\ \text{is strictly increasing}\},\\
\mathcal{K}_\infty&:=\{\gamma\in \mathcal{K}\mid\gamma\  \text{is unbounded}\},\\
\mathcal{L}&:=\{\gamma:\mathbb{R}_+\rightarrow\mathbb{R}_+\mid\gamma\ \text{is continuous, strictly}\\
&\text{\qquad \qquad \qquad \qquad \quad decreasing with}\ \lim_{r\rightarrow\infty}\gamma(r)=0 \},\\
\mathcal{KL}&:=\{\beta: \mathbb{R}_+\times \mathbb{R}_+\rightarrow\mathbb{R}_+\mid\beta\ \text{is continuous}, \beta(\cdot,t)\in\mathcal{K},\\
&\qquad \qquad \qquad \qquad \quad\quad \quad \ \forall t\geq0, \beta(r,\cdot)\in\mathcal{L}, \forall r>0\}.
\end{align*}

\section{Problem setting}\label{sec:pro setting}

In this paper, we use the following notion of a dynamical system with outputs.

\begin{definition}\label{dynamical system}
A \emph{dynamical system} $\Sigma=(X, Y,\phi,h)$ with outputs consists of
\begin{enumerate}[(i)]
\item  a normed linear space $(X, \|\cdot\|_X)$, called the \emph{state space};

\item  a normed linear space  $(Y, \|\cdot\|_Y)$, called the \emph{output value space};

\item  a continuous map $h: X \rightarrow {Y}$, called the \emph{output map};

\item  a map $\phi: D_\phi\rightarrow X$ with $D_\phi\subseteq\mathbb{R}_+\times X$, called the \emph{transition map}.

For every $x\in X$, there exists $t_m(x)\in(0,\infty]$, called the \emph{maximal existence time} of the trajectory starting from $x$, such that
\begin{equation*}D_\phi=\bigcup_{x\in X}[0,t_m(x))\times \{x\}.\end{equation*}

Furthermore, the following properties hold:
   \begin{enumerate}
\item [($\Sigma$1)] \emph{The identity property}: for every $x\in X$ it holds that $\phi(0,x)=x$.
  \item [($\Sigma$2)] \emph{Continuity}: for every $x\in X$ the map $t\mapsto\phi(t,x)$ is continuous on $[0,t_m(x))$.
      \item [($\Sigma$3)] \emph{The semigroup property}: for every $x\in X$, for every $t\in[0,t_m(x))$ and every $s\ge0$ such that $t+s<t_m(x)$, one has 
         $s<t_m(\phi(t,x))$ and $\phi(s,\phi(t,x))=\phi(t+s,x)$.
\end{enumerate}
\end{enumerate}

 For each $x\in X$, the corresponding output trajectory is defined by
\[
y(t,x):=h(\phi(t,x)),\quad t\in[0,t_m(x)).
\]

\end{definition}

\begin{definition}
    A dynamical system $\Sigma=(X,Y,\phi,h)$ is called \emph{forward complete} if $t_m(x)=\infty$ for every $x\in X$.
\end{definition}

\begin{remark}
    We emphasize that the state space $X$ is only assumed to be a normed linear space and need not be complete. This generality is important for the PDE example in Section~\ref{sec:example}, where a non-complete state space is used.
\end{remark}

In the following, unless stated otherwise, $\Sigma$ denotes a dynamical system
with outputs.
Since the stability notions below are formulated with respect to the origin, we
assume that $0$ is an \emph{equilibrium} of $\Sigma$, i.e.,
\begin{equation}\label{eq:equilibrium}
t_m(0)=\infty,\quad \phi(t,0)=0,\quad t\ge0.
\end{equation}

The abstract framework encompasses, in particular, suitable classes of autonomous ordinary differential equations, time-delay systems, and many classes of partial differential equations and  evolution equations on Banach spaces.

\begin{definition}\label{def:BIC}
    We say that $\Sigma$ satisfies the \emph{boundedness-implies-continuation (BIC)} property, if for every $x\in X$
    \begin{align*}
        \sup_{t\in[0,t_m(x))}\|\phi(t,x)\|_X<\infty \quad \Rightarrow \quad t_m(x)=\infty.
    \end{align*}
\end{definition}

We proceed to the main concept of this work.

\begin{definition}\label{def:abOSS}
  System  $\Sigma$ is   called \emph{output-to-state stable (OSS)}, if there exist $\beta\in\mathcal{KL}$ and $\gamma\in\mathcal{K}$ such that for all $ x\in X$ and all $t\in[0,t_m(x))$ it holds that
\begin{align}
\label{abOSS}
\|\phi(t,x)\|_{X}\leq \beta(\|x\|_{X},t)+ \gamma(\max_{s\in[0,t]}\|y(s,x)\|_Y).
\end{align}
\end{definition}

The OSS property not only ensures that bounded outputs imply bounded states on the maximal interval of existence, but also guarantees that vanishing outputs lead to vanishing states under the BIC property, as shown  in Corollary~\ref{cor:OSS-has-vanishing-output-vanishing-state}.
This plays a vital role in observer theory and detectability analysis.

We use OSS as an abbreviation for both \q{output-to-state stability} and \q{output-to-state stable}.

\begin{remark}
\label{remark: output-oss}
 For a dynamical system $\Sigma$,  each of the following conditions is sufficient for  OSS:
 \begin{enumerate}
 \item [(i)] If  there exists $\beta\in\mathcal{KL}$ such that 
   \[\|\phi(t,x)\|_X\leq \beta(\|x\|_X,t),\quad x \in X, \quad t\in[0,t_m(x)),\] 
   then $\Sigma$ is OSS with any $\gamma \in \mathcal{K}$.

 \item[(ii)]   Assume that there exists $\xi\in\mathcal{K}_\infty$ such that
\[
\|h(x)\|_Y \geq \xi(\|x\|_X), \quad x\in X.
\]
Then $\Sigma$ is OSS with $\gamma=\xi^{-1}$.

 \end{enumerate}
\end{remark}

\section{Direct OSS Lyapunov theorem}\label{sec:direct OSS}

In this section, we introduce coercive OSS Lyapunov functions in implication and dissipation formulations, study the relation between them, and establish a direct Lyapunov theorem for OSS. We begin by recalling the concept of the upper right-hand Dini derivative.

For a function $q:\mathbb{R}\to\mathbb{R}$, define the (upper right-hand) \emph{Dini derivative} of $q$ at $t$ by
\begin{align*}
D^+q(t)
=\mathop{\overline{\lim}}\limits_{\tau\rightarrow+0}
\frac{q(t+\tau)-q(t)}{\tau}.
\end{align*}
 The limsup is understood in the extended real line, so $D^+q(t)$ may equal $+\infty$ or $-\infty$.
Given $V:X\to \R_+$, define the \emph{Lie derivative} of $V$ at $x \in X$ by
\[
\dot{V}(x):=D^+V(\phi(t,x))|_{t=0}=
\mathop{\overline{\lim}}\limits_{\tau\rightarrow+0}
\frac{V(\phi(\tau,x))-V(x)}{\tau},
\]
    provided that $t\mapsto V(\phi(t,x))$ is continuous at $t=0$.

\begin{definition}
\label{def:abOSS LF-dissipative}
A   function $V: X\rightarrow \mathbb{R}_{+}$ is called a  \emph{coercive OSS Lyapunov function in dissipation form} for  $\Sigma$ if:

 \begin{enumerate}
\item[(i)] There exist $\psi_1, \psi_2\in \mathcal{K}_{\infty}$  such that for every $x\in X$:
\begin{align}\label{eq:aboss LF1}
\psi_1(\|x\|_{X})\leq V(x)\leq \psi_2(\|x\|_{X}).
\end{align}
\item[(ii)] For every $x\in X$,  the map $t\mapsto V(\phi(t,x))$, $t\in[0,t_m(x))$, is continuous.
\item[(iii)] There exist $\alpha, \sigma\in\mathcal{K}_{\infty}$ such that for every $x\in X$,
the following dissipation inequality holds:
\begin{align}\label{aboss LF2}
\dot{V}(x)\leq-\alpha(\|x\|_{X})+\sigma(\|h(x)\|_{Y}).
\end{align}
\end{enumerate}

The function $V$ is called a  \emph{coercive OSS Lyapunov function in  implication form} for $\Sigma$, if
items (i) and (ii)  hold and
\begin{enumerate}
    \item[(iii')] there exist $\tilde{\alpha}\in \mathcal{P}$ and $\chi\in\mathcal{K}_{\infty}$ such that for every $x\in {X}$, the following implication holds:
\begin{equation}
\label{eq:aboss Implication}
\hspace{-4mm}\|x\|_{X} \geq\chi(\|h(x)\|_{Y})
\quad \Rightarrow\quad  \dot{V}(x)\leq-\tilde{\alpha}(\|x\|_{X}).
\end{equation}
\end{enumerate}

\end{definition}

Next we clarify the relation between the  dissipation and implication formulations of OSS Lyapunov functions.
\begin{proposition}
\label{prop:dflf-criterion}
A function
$V$ is a coercive  OSS Lyapunov function in dissipation form for $\Sigma$
if and only if there is $\chi\in\Kinf$ such that  two conditions hold:
\begin{itemize}
    \item[(i)] $V$ is a coercive OSS Lyapunov function in implication form for $\Sigma$ with $\chi$ and $\tilde{\alpha}\in \mathcal{K}_\infty$.
    \item[(ii)] There is $ \alpha_1\in\mathcal{K}_\infty$ so that 
    for every $x\in X$, we have:
    \begin{align}\label{eq:implication2}
    \hspace{-2mm}
        \|x\|_X < \chi(\|h(x)\|_Y)\srs \dot{V}(x)\leq \alpha_1(\|h(x)\|_Y).
    \end{align}  
\end{itemize}
\end{proposition}

\begin{proof}
\q{$\Rightarrow$}. Let $V$ be a coercive OSS Lyapunov function in dissipation form with $\psi_1,\psi_2,\alpha,\sigma\in \mathcal{K}_\infty$. Define $\chi(r):=\alpha^{-1}(2\sigma(r))$,  $r\in\mathbb{R}_+$.
Clearly, (ii) holds with $\alpha_1=\sigma$. To show (i), note that
\begin{equation*}
\|x\|_{X} \geq\chi(\|h(x)\|_{Y})
 \sis \sigma(\|h(x)\|_{Y})\leq \frac{1}{2}\alpha(\|x\|_{X}).
\end{equation*}
By \eqref{aboss LF2}, we have
\begin{equation*}
 \sigma(\|h(x)\|_{Y})\leq \frac{1}{2}\alpha(\|x\|_{X})
 \qrq \dot{V}(x)\leq-\frac{1}{2}\alpha(\|x\|_{X}).
 \end{equation*}
Therefore, $V$ is a coercive OSS Lyapunov function in implication form with $\tilde{\alpha}=\frac{1}{2}\alpha$.

\q{$\Leftarrow$}. 
Let $V$ be a coercive OSS Lyapunov function in implication form with $\psi_1, \psi_2, \chi,\tilde{\alpha}\in \mathcal{K}_\infty$. We have
 \begin{align*}
        \|x\|_X\geq \chi(\|h(x)\|_Y) \qrq \dot{V}(x)\leq -\tilde{\alpha}(\|x\|_X).
    \end{align*}
Furthermore, assume that \eqref{eq:implication2} holds for $\alpha_1\in\mathcal{K}_\infty$. That is, 
\begin{align*}
       \|x\|_X&<\chi(\|h(x)\|_Y)\\
        &\Rightarrow\dot{V}(x)\leq -\|x\|_X + \|x\|_X + \alpha_1(\|h(x)\|_Y)\\
        &\qquad\quad~~ \leq -\|x\|_X + \chi(\|h(x)\|_Y) + \alpha_1(\|h(x)\|_Y).
\end{align*}

Combining the two cases, we obtain the dissipation inequality for every $x\in X$, with $\alpha(r)=\min\{\tilde{\alpha}(r),r\}$ and $\sigma=\alpha_1 + \chi$.
Thus, $V$ is a coercive OSS Lyapunov function in dissipation form. 
\end{proof}

\begin{remark}
 In Proposition~\ref{prop:dflf-criterion}, the additional condition (ii) is similar to the bound used in the definition of strong implication-form ISS-Lyapunov functions for discontinuous discrete-time systems \cite [Definition 4.1]{grune2014iss}. The equivalence between the existence of this type of Lyapunov functions and Lyapunov functions in dissipation form is analogous to Proposition \ref{prop:dflf-criterion}; see \cite [Theorem 4.4(a)]{grune2014iss}.
\end{remark}

\begin{remark}
If there exists $\alpha_2\in\mathcal{K}_\infty$ such that $\dot{V}(x)\le \alpha_2(\|x\|_X)$ for all $x\in X$,
then condition (ii) in Proposition~\ref{prop:dflf-criterion} is satisfied with $\alpha_1=\alpha_2\circ\chi$. Consequently, $V$ is an OSS Lyapunov function in dissipation form if and only if it is an OSS Lyapunov function in  implication form with $\tilde{\alpha}\in\mathcal{K}_\infty$. In particular, this is the case for ODE systems with  locally Lipschitz right-hand sides  and $C^1$ functions $V$. This is analogous to the equivalence between the implication and dissipation formulations of ISS Lyapunov functions in the finite-dimensional setting; see \cite[Proposition 2.21]{mironchenko2023}.
\end{remark}

\begin{remark}\label{re:im-not-di}
    The counterexample in \cite[Proposition~5.2]{chen2026detectability} shows that, even for infinite-dimensional  linear systems with a bounded output operator, the existence of a coercive OSS Lyapunov function in implication form does not, in general, imply the existence of one in dissipation form.
\end{remark}

Proposition~\ref{prop:dflf-criterion} establishes that, under suitable conditions, the dissipation form can be derived from an OSS Lyapunov function in implication form with a decay rate in \( \mathcal{K}_\infty \).

The following proposition considers the case where the decay rate belongs only to \( \mathcal{P} \), and shows that coercivity enables a transformation that lifts the decay rate into \( \mathcal{K}_\infty \). 
This proposition is a direct counterpart of \cite[Remark 4.1]{LSW96}; see also  \cite[Proposition 2.17]{mironchenko2023}.

\begin{proposition}\label{prop:p to k}
Let \(V\) be a coercive OSS Lyapunov function in implication form with a decay rate $\tilde{\alpha}\in\mathcal{P}$. Then there exists  \( \xi \in \mathcal{K}_\infty \cap C^1(\mathbb{R}_+, \mathbb{R}_+) \), smooth on \( (0, \infty) \),  such that $\xi \circ V$ is a coercive OSS Lyapunov function in implication form with a decay rate belonging to
$\mathcal{K}_\infty$.
\end{proposition}
\begin{proof}
    Let \(V\) be  a coercive OSS Lyapunov function in implication form with  $\psi_1,\psi_2,\chi\in\mathcal{K}_\infty$ and 
 a decay rate $\tilde{\alpha}\in\mathcal{P}$.   Note that $\tilde{\alpha}\circ\psi_1^{-1}\in \mathcal{P}$. 
 By \cite[Proposition~A.13]{mironchenko2023}, there exist $k_1\in\mathcal{K}_\infty$ and $l_1\in\mathcal{L}$ such that
\begin{align*}
    \tilde{\alpha}(\psi_1^{-1}(r))\geq k_1(r)l_1(r),\quad r\in\mathbb{R}_+.
\end{align*}
Without loss of generality, assume that $l_1$ is smooth on 
\((0, \infty)\); see \cite[Lemma 2]{Kel14}. Now define 
\begin{align*}
    \xi(r):=\int_0^r\frac{1}{l_1(s)}\, \mathrm{d}s,\quad r
    \in \mathbb{R}_+,
\end{align*}
and $W:=\xi\circ V$. Since $\xi$ is continuous and $t\mapsto V(\phi(t,x))$ is continuous for every $x\in X$, the map $t\mapsto W(\phi(t,x)) $ is continuous on $[0,t_m(x))$. For $x\in X$, we deduce
\begin{align*}
\xi \circ \psi_1(\|x\|_X) \leq W(x) \leq \xi \circ \psi_2(\|x\|_X).
\end{align*}

For $x\neq0$, apply  Lemma~\ref{le:chain-rule} to get $\dot{W}(x)=\frac{1}{l_1(V(x))}\dot{V}(x)$. At $x=0$, both Lie derivatives vanish because zero is an equilibrium. Thus, 
\begin{align*}
\|x\|_X\geq \chi(\|h(x)\|_Y)\Rightarrow \dot{W}(x)
\leq&-\frac{1}{l_1(\psi_1(\|x\|_X))}\tilde{\alpha}(\|x\|_X)\\
\leq&-k_1\circ \psi_1(\|x\|_X),
\end{align*}
and the proof is complete.
\end{proof}

We are  ready to state the direct Lyapunov theorem for OSS.

 \begin{theorem}[Direct coercive OSS Lyapunov theorem]\label{thm:abdirect OSS LT}
     If there exists a coercive OSS Lyapunov function (in dissipation or implication form) for $\Sigma$, then $\Sigma$ is OSS.
 \end{theorem}

\begin{proof}
By Proposition~\ref{prop:dflf-criterion}, it suffices  to consider coercive OSS Lyapunov functions in implication form. Let $\Sigma$ admit a coercive OSS Lyapunov function in implication form with $\psi_1,\psi_2,\chi\in\mathcal{K}_\infty$ and $\tilde{\alpha}\in\mathcal{P}$. In view of Proposition~\ref{prop:p to k}, we assume without loss of generality that $\tilde{\alpha} \in \mathcal{K}_\infty$. 

Fix $x\in X$ and consider $\phi(\cdot):=\phi(\cdot,x)$.
Define
\begin{align*}
a(t):=\psi_2\circ\chi(\max_{s\in[0,t]}\|y(s,x)\|_Y),\quad t\in[0,t_m(x))
\end{align*}
and 
\begin{align*}
    S(t):=\{z\in X \mid V(z)\leq a(t)\}, \quad t\in[0,t_m(x)).
\end{align*}

\textbf{Step 1.} We are going to prove that, for every $T\in[0,t_m(x))$,
\begin{align*}
     \phi(T)\in S(T) \qrq \phi(t)\in S(t),\quad t\in[T,t_m(x)).
\end{align*}

Assume $T\in[0,t_m(x))$ is such that $\phi(T)$ lies in $S(T)$. To obtain a contradiction, suppose that there exist some $\varepsilon>0$ and $t_1\in [T,t_m(x))$  such that $V(\phi(t_1))>a(t_1)+\varepsilon$. Define
\begin{align*}
    T_1:=\inf\{t\in[T,t_m(x))\mid V(\phi(t))\geq a(t)+\varepsilon\}.
\end{align*}
By continuity, $T<T_1\le t_1$ and $V(\phi(T_1))=a (T_1)+\varepsilon$.
Since $t\mapsto \phi(t,x)$ is continuous on $[0,t_m(x))$ and $h$ is continuous, $t\mapsto y(t,x)$ is continuous on $[0,t_m(x))$. Hence $t\mapsto \max_{s\in[0,t]}\|y(s,x)\|_Y$ is continuous on $[0,t_m(x))$, and therefore $ a(\cdot)
$ is continuous on this interval. In addition, the map $t\mapsto V(\phi(t))$ is continuous on $[0,t_m(x))$. Therefore, there is a neighborhood  $\mathcal{N}(T_1)\subset[0,t_m(x))$ of $T_1$ such that $V(\phi(t))\geq a(t)$ for all $t\in \mathcal{N}(T_1)$. From \eqref{eq:aboss LF1}, we have $\psi_2(\|\phi(t)\|_X)\geq a(t)$ for  $t\in \mathcal{N}(T_1)$. By the definition of $a$,  
\begin{align*}
\|\phi(t)\|_X\geq \chi(\max_{s\in[0,t]}\|y(s,x)\|_Y) \ge \chi(\|y(t,x)\|_Y) 
\end{align*} holds for $t\in \mathcal{N}(T_1)$. By \eqref{eq:aboss Implication}, we also have $\dot V(\phi(t))\leq-\tilde{\alpha}(\|\phi(t)\|_{X})$, $t\in \mathcal{N}(T_1).$

By \cite[Theorem~XI.4.1]{Bruckner1978},  the negativity of Dini derivative guarantees that $t\mapsto V(\phi(t))$ is nonincreasing on $\mathcal{N}(T_1)$. Observe that $a(\cdot)$ is nondecreasing on $ [0,t_m(x))$. Hence, 
\begin{align*}
V(\phi(t))\geq V(\phi(T_1))\geq a(T_1)+\varepsilon \geq a(t)+\varepsilon
\end{align*}
for all $t\in \mathcal{N}(T_1)\cap (T, T_1) $, contradicting the choice of $T_1$.

Let $T^*:=\inf \{t\in[0,t_m(x))\mid  \phi(t)\in S(t)\}$,  with the convention that $T^*:=t_m(x)$ if the set is empty. If $T^*<t_m(x)$, continuity implies that $\phi(T^*)\in S(T^*)$. Thus, by Step 1, $\phi(t)\in S(t)$ for all $t\in[T^*,t_m(x)).$  It follows from  \eqref{eq:aboss LF1} that for all 
$t\in[T^*,t_m(x))$
\begin{align}\label{eq:abstep1}
    \|\phi(t)\|_X\leq\psi_1^{-1}\circ\psi_2 \circ\chi(\max_{s\in[0,t]}\|y(s,x)\|_Y).
\end{align}

\textbf{Step 2.}  By the definition of $T^*$, \eqref{eq:aboss LF1} and \eqref{eq:aboss Implication}, it follows that
\begin{equation*}
\dot V(\phi(t))
\leq- \tilde{\alpha}\circ\psi_2^{-1}(V(\phi(t))), \quad t \in [0, T^*).
\end{equation*}
According to  \cite[Proposition~A.35]{mironchenko2023}, there exists $\beta_0\in\mathcal{KL}$,  independent of the initial state $x$, such that
\begin{equation*}
V(\phi(t))\leq \beta_0(V(\phi(0)),t),\quad t \in[0,T^*).
\end{equation*}
It follows from  \eqref{eq:aboss LF1} that
\begin{equation}\label{eq:abstep2}
\|\phi(t)\|_X\leq\psi_1^{-1}\circ\beta_0(\psi_2(\|{x}\|_X),t),\quad t \in[0,T^*).
\end{equation}
Combining \eqref{eq:abstep1} and \eqref{eq:abstep2}, we obtain \eqref{abOSS} with $\beta(r,t)=\psi_1^{-1}\circ\beta_0(\psi_2(r),t)$ and $\gamma(r)=\psi_1^{-1}\circ\psi_2(\chi(r))$. Since $x\in X$ was arbitrary and $\beta$ and $\gamma$ are independent of $x$, the estimate holds  for all $x\in X$. Thus,  $\Sigma$ is OSS.
\end{proof}

\section{Uniform global asymptotic stability modulo output}\label{sec:gasmo}

In this section, we recall the following notion from \cite[Definition 2.9]{krichman2001input} and prove its equivalence to OSS.

\begin{definition}\label{UGASMO}
System $\Sigma$ is called  \emph{uniformly globally asymptotically stable modulo output (UGASMO)}, if there exist  $\rho\in \mathcal{K}_{\infty} $ and $\theta\in \mathcal{KL}$
such that for all $x\in X$ and all $T\in[0,t_m(x))$, if
\begin{align*}
\|\phi(t,x)\|_X\geq\rho(\|y(t,x)\|_Y)\quad\forall t\in[0,T],
\end{align*}
then
\begin{align}\label{eq:theta}
  \|\phi(t,x)\|_X\leq\theta(\|x\|_X,t) \quad\forall t\in[0,T].
\end{align}
\end{definition}

\begin{figure}[htbp] 
\centering 
\begin{tikzpicture}[ 
    x=0.75cm, 
    y=1.0cm, 
    >=Stealth, 
    font=\small 
]


\def\tstar{5.15} 
\def\xend{10.1} 
\def\ymax{5.15}


\fill[green!8] 
    (0,0) rectangle (\tstar,\ymax); 
 
\fill[red!8] 
    (\tstar,0) rectangle (\xend,\ymax); 


\draw[->, thick] 
    (0,0) -- (10.55,0) 
    node[right, font=\large] {$t$}; 
 
\draw[->, thick] 
    (0,0) -- (0,5.25); 
 
\node[below left] at (0,0) {$0$};


\draw[dashed, thick] 
    (\tstar,0) -- (\tstar,\ymax); 
 
\draw (\tstar,0.06) -- (\tstar,-0.06); 
 
\node[below=5pt, font=\large] 
    at (\tstar,0) {$t^{*}$};


\draw[ 
    green!55!black, 
    very thick, 
    dash pattern=on 8pt off 3pt on 1.5pt off 3pt 
] 
plot[smooth] coordinates { 
    (0.00,4.72) 
    (0.45,4.10) 
    (0.95,3.55) 
    (1.55,3.05) 
    (2.15,2.65) 
    (2.85,2.25) 
    (3.55,1.92) 
    (4.25,1.62) 
    (5.15,1.28) 
    (6.00,1.03) 
    (6.90,0.83) 
    (7.80,0.67) 
    (8.65,0.54) 
    (10.1,0.46) 
};

\draw[red, very thick]
plot[smooth] coordinates {
    (0.00,3.62)
    (0.35,3.28)
    (0.70,2.82)
    (1.05,2.35)
    (1.45,1.88)
    (1.85,1.42)
    (2.25,1.58)
    (2.70,1.56)
    (3.05,1.57)
    (3.45,1.47)
    (3.85,1.26)
    (4.25,1.07)
    (4.65,0.94)
    (5.15,0.82)
    (5.55,0.62)
    (6.00,0.48)
    (6.45,0.68)
    (6.95,0.82)
    (7.35,0.42)
    (7.85,0.34)
    (8.30,0.86)
    (8.75,0.74)
    (9.20,0.28)
    (9.65,0.24)
    (10.10,0.58)
};

\draw[
    blue,
    very thick,
    dashed
]
plot[smooth] coordinates {
    (0.00,1.68)
    (0.35,1.70)
    (0.75,1.65)
    (1.15,1.52)
    (1.60,1.34)
    (2.10,1.17)
    (2.65,1.02)
    (3.25,0.91)
    (3.85,0.85)
    (4.45,0.83)
    (5.15,0.82)
    (5.55,0.96)
    (6.00,1.0)
    (6.45,1.52)
    (6.95,1.74)
    (7.35,2.15)
    (7.85,2.2)
    (8.30,2.15)
    (8.75,1.98)
    (9.20,1.12)
    (9.65,0.91)
    (10.10,0.76)
};

\node[ 
    green!55!black, 
    font=\large 
] 
at (2.9,3.0) 
{\small$\theta(\|x\|_{X},t)$};

\node[ 
    red, 
    font=\large 
] 
at (2.6,1.8) 
{\small$\|\phi(t,x)\|_{X}$};

\node[ 
    blue, 
    font=\large 
] 
at (2,0.65) 
{\small $\rho(\|y(t,x)\|_{Y})$};

\node[ 
    font=\bfseries\large, 
    align=center 
] 
at (2.7,4.7) 
{\scriptsize State-dominated region\\[1pt]
\footnotesize$\|\phi(t,x)\|_{X} \ge \rho(\|y(t,x)\|_{Y})$};

\node[ 
    font=\bfseries\large, 
    align=center 
] 
at (7.7,4.7) 
{\scriptsize Output-dominated region\\[1pt]
\footnotesize$\|\phi(t,x)\|_{X} < \rho(\|y(t,x)\|_{Y})$};

\end{tikzpicture} 
 
\caption{Illustration of the UGASMO property 
} 
\label{fig:ugasmo-overview} 
\end{figure}
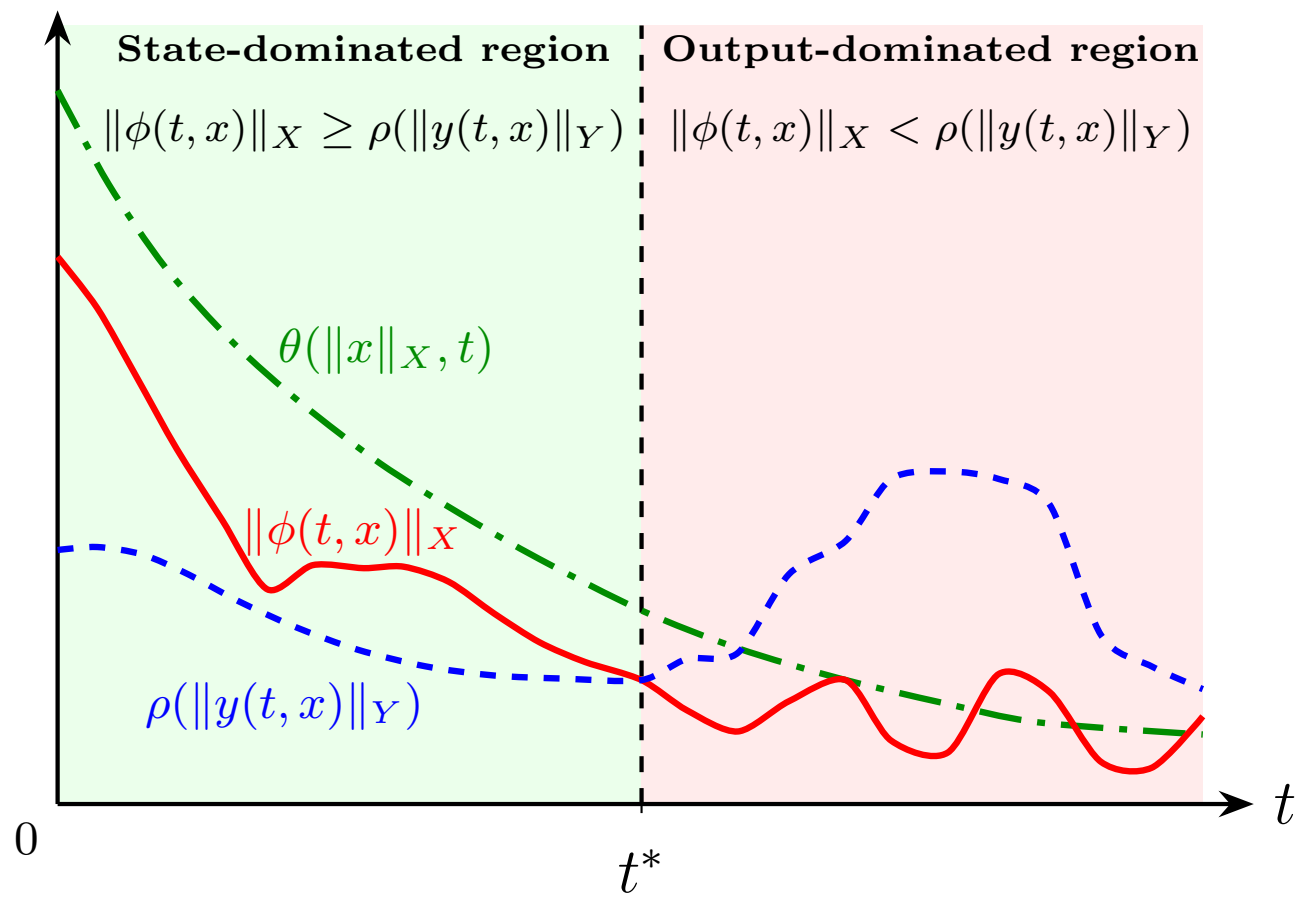

As illustrated in Figure~\ref{fig:ugasmo-overview}, the state is required to satisfy a uniform $\mathcal{KL}$ decay estimate as long as the trajectory remains in the state-dominated region. Once the trajectory enters the output-dominated region, the premise of the UGASMO estimate is no longer satisfied, and the $\mathcal{KL}$ estimate can no longer be concluded.

 An equivalent criterion to UGASMO was established in  \cite[Proposition 2.10]{krichman2001input} for ODE systems. Here, we extend it to the infinite-dimensional setting.

 \begin{proposition}\label{prop iff}
System $\Sigma$ is UGASMO if and only if there exists $\rho\in\mathcal{K}_{\infty}$ so that the following statements hold:
\begin{enumerate}
\item[(i)] There exists $\kappa\in \mathcal{K}$ such that for all $x\in X$ and all  $T\in[0,t_m(x))$ such that
\begin{align*}
\|\phi(t,x)\|_X\geq\rho(\|y(t,x)\|_Y)\quad\forall t\in[0,T],
\end{align*}
the following estimate holds:
\begin{align*}
  \|\phi(t,x)\|_X\leq\kappa(\|x\|_X) \quad\forall t\in[0,T].
\end{align*}

\item[(ii)] For all $\varepsilon>0$, all $r>0$, there exists  $T(r,\varepsilon)>0$ such that for all $x\in X$ with $\|x\|_X\leq r$ and all $T\in[T(r,\varepsilon),t_m(x))$, if
\begin{align*}
\|\phi(t,x)\|_X\geq\rho(\|y(t,x)\|_Y)\quad\forall t\in[0,T],
\end{align*}
 then
 \begin{align*}
 \|\phi(t,x)\|_X<\varepsilon \quad \forall t\in[T(r,\varepsilon),T].
 \end{align*}
\end{enumerate}
 \end{proposition}

\begin{proof}
The necessity is obvious. To prove sufficiency, assume that statements (i) and (ii) hold. Define the function \( g: \mathbb{R}_+ \times \mathbb{R}_+ \to \mathbb{R} \) by
\begin{align*}
g(s,t):=\sup\{&\|\phi(t,x)\|_X\mid \|x\|_X\leq s, t<t_m(x),\\
&\|\phi(\hat{t},x)\|_X\geq\rho(\|y(\hat{t},x)\|_Y)\quad\forall \hat{t}\in[0,t] \},
\end{align*}
where, by convention, the supremum of the empty set is defined to be zero. By (i), $g(s,t)\le \kappa(s)$ for all $s,t\in\mathbb{R}_+$, and (ii) implies that, for every $r,\varepsilon>0$,
$g(s,t)\leq\varepsilon$
whenever $s\leq r$ and $t\geq T(r,\varepsilon)$. Since $\kappa\in\mathcal{K}$, it can be majorized by a function in $\mathcal{K}_\infty.$  Thus all the assumptions of  \cite[Proposition 4.5]{Mir25b} are satisfied. Then there is $\beta\in \mathcal{KL}$ such that
$g(s,t) \leq \beta(s,t)$ for all $ s,t \in\mathbb{R}_+.$ By the definition of $g$, this yields
$\|\phi(t,x)\|_X\le g(\|x\|_X,t) \le \beta(\|x\|_X,t)$ for all $x\in X$ and $t\in[0,t_m(x))$ satisfying the premise in statement (i).
\end{proof}

\begin{remark}
For forward complete systems, if the outputs are identically zero, then the statement (i) of Proposition \ref{prop iff} becomes uniform global stability and the statement (ii) of Proposition \ref{prop iff} is uniform global attractivity; see \cite[Definition B.21]{mironchenko2023}.
\end{remark}

Now we can establish the equivalence between OSS and UGASMO. The proof of the implication from OSS to UGASMO follows the approach of \cite[Lemma 2.3]{sontag1996}, while the converse implication is established directly.

\begin{theorem}
    \label{thm:oss-equi-ugasmo}
 System $\Sigma$ is OSS $\iff$ $\Sigma$ is UGASMO.

\end{theorem}

\begin{proof}
($\Rightarrow$)
Suppose that $\Sigma$ is  OSS with $\beta$ and $\gamma$. We prove the claim by using the characterization of UGASMO obtained in Proposition \ref{prop iff}.
Let $\hat{\kappa}(s):=\beta(s,0)$, $s \geq 0$. Without loss of generality, assume that $\gamma\in \mathcal{K}_\infty$ in \eqref{abOSS} and $\hat{\kappa}(s)>s$ for all $s>0$.  Pick $\rho\in\mathcal{K}_\infty$  satisfying
\begin{align}
\label{eq:choice-rho}
\rho(s)>\hat{\kappa}(4\gamma(s))+\gamma(s),\quad s>0.
\end{align}

{By applying $\rho^{-1}$ to} both sides and substituting $s=\gamma^{-1}(\frac{r}{4})$, we see that \eqref{eq:choice-rho} is equivalent to
\begin{align}
\label{eq:choice-rho-2}
\gamma\circ \rho^{-1}\bigl(\hat{\kappa}(r)+\frac r4\bigr)<\frac r4,\quad r>0.
\end{align}

\textbf{Step 1}.  We prove that the condition (i) of Proposition \ref{prop iff} holds.
  Pick any $ x \in X$ and any  $T\in [0,t_m(x))$ with $\|\phi(t,x)\|_X\geq\rho(\|y(t,x)\|_Y)$ for all $0\leq t\leq T$.
The case $x=0$ is immediate, so assume that $x\neq0.$

Note that $\|x\|_X\le \hat{\kappa}(\|x\|_X)+\frac{\|x\|_X}{4}.$ Then by the definition of $\hat{\kappa}$, for $t=0$ we have
\begin{align*}
\gamma(\|y(0,x)\|_Y)&\leq\gamma(\rho^{-1}(\|x\|_X))\\
&\leq\gamma\bigl(\rho^{-1}\bigl(\hat{\kappa}(\|x\|_X)+\frac {\|x\|_X}{4}\bigr)\bigr)\overset{\eqref{eq:choice-rho-2}}{<}\frac{\|x\|_X}{4}.
\end{align*}

Let $t_1:=\inf\{t\in(0,T] \mid \gamma(\|y(t,x)\|_{Y})\geq\frac{\|x\|_X}{4}\}$. By the previous argument, $t_1>0$.  As $y$ is continuous in time, we have $\gamma(\|y(t_1,x)\|_{Y})=\frac{\|x\|_X}{4}$ and
\begin{equation*}
 \gamma(\|y(t,x)\|_{Y})<\frac{\|x\|_X}{4}\quad \forall t\in[0,t_1).
\end{equation*}
By the definition of OSS and since $\hat{\kappa}(s)>s$ for all $s>0$, we also get for all $t\in[0,t_1]$
\begin{align}\label{eq:hat_kapaa}
\|\phi(t,x)\|_X&\leq\beta(\|x\|_{X},t)+\gamma(\max_{s\in[0,t]}\|y(s,x)\|_Y)\notag\\ 
&\le\hat{\kappa}(\|x\|_X)+\frac{\|x\|_X}{4}.
\end{align}

Using \eqref{eq:choice-rho-2} again, we obtain for all $t\in[0,t_1]$
\begin{align*}
\gamma(\|y(t,x)\|_Y)&\leq\gamma(\rho^{-1}(\|\phi(t,x)\|_X))\\
&\overset{\eqref{eq:hat_kapaa}}{\leq}\gamma\bigl(\rho^{-1}\bigl(\hat{\kappa}(\|x\|_X)+\frac {\|x\|_X}{4}\bigr)\bigr)\overset{\eqref{eq:choice-rho-2}}{<}\frac{\|x\|_X}{4},
\end{align*}
which contradicts the definition of $t_1$. This contradiction shows that the set defining $t_1$ is empty. Thus,
\begin{equation*}\gamma(\|y(t,x)\|_Y)<\frac{\|x\|_X}{4}, \quad t\in[0,T],\end{equation*}
and hence 
\begin{align}\label{eq:gamma-4}
\gamma(\max_{s\in[0,t]}\|y(s,x)\|_Y)\le \frac{\|x\|_X}{4}, \quad t\in[0,T].
\end{align}

Therefore, by \eqref{eq:hat_kapaa}, we obtain
\begin{align*}
    \|\phi(t,x)\|_X\le \hat{\kappa}(\|x\|_X)+\frac{\|x\|_X}{4},\quad t\in [0,T].
\end{align*}
Thus condition (i) of Proposition \ref{prop iff} holds with $\kappa(r)=\hat{\kappa}(r)+\frac r4$.

\textbf{Step 2}. We prove that condition (ii) of Proposition \ref{prop iff} holds.
Fix $\varepsilon>0$ and $r>0$. Let $k=k(r,\varepsilon)$ be a positive integer satisfying $2^{-k}r<\varepsilon$  and define $r_i:=2^{1-i}r$ for $i=1,2,\ldots, k+1$. Then
\begin{equation*}
r_{i+1}=\frac{r_i}{2},\quad i=1,\cdots,k\quad\text{and}\quad r_{k+1}<\varepsilon.
\end{equation*}

For each $i=1,...,k,$ choose $T_{r_i}>0$  such that
\begin{align}\label{eq:beta-ri}
    \beta(r_i,t)<\frac{r_i}{4},\quad t\geq T_{r_i}.
\end{align}
Now define \begin{equation*}T(r,\varepsilon) :=T_{r_1}+T_{r_2}+\cdots+T_{r_k}.\end{equation*}

Let $x\in X$ satisfy $\|x\|_X\leq r$, and let $T\in [0, t_m(x))$ satisfy $T\geq T(r,\varepsilon)$. Assume that
\begin{equation*}\|\phi(t,x)\|_X\geq\rho(\|y(t,x)\|_Y),\quad t\in[0,T].\end{equation*}

In Step 1, we have shown that this implies \eqref{eq:gamma-4} and  so
\begin{equation*}\gamma(\max_{s\in[0,t]}\|y(s,x)\|_Y)\le\frac{\|x\|_X}{4}\leq\frac{r_1}{4},\quad t\in[0,T].\end{equation*}
 
Then, for every $t\in[T_{r_1}, T]$, we have from \eqref{abOSS}
\begin{align*}
\|\phi(t,x)\|_X&\leq\beta(\|x\|_X,t)+\gamma(\max_{s\in[0,t]}\|y(s,x)\|_Y)\\
&\leq\beta(r_1,t)+\gamma(\max_{s\in[0,t]}\|y(s,x)\|_Y)\\
&\overset{\eqref{eq:beta-ri}}{\leq} \frac{r_1}{4}+\frac{r_1}{4}=r_2.
\end{align*}

By the semigroup property, the trajectory starting from $\phi(T_{r_1},x)$ satisfies the assumption of Step 1 on $[0,T-T_{r_1}]$. Hence, applying \eqref{eq:gamma-4} again, we obtain for all $t\in[0,T-T_{r_1}]$
\begin{equation*}\gamma\big(\max_{s\in[0,t]}\|y(s,\phi(T_{r_1},x))\|_Y\big)\le\frac{\|\phi(T_{r_1},x)\|_X}{4}\leq\frac{r_2}{4}.\end{equation*}

Thus,  for every $t\in[T_{r_2}, T-T_{r_1}]$, we have
\begin{align*}
\|\phi(t,\phi(T_{r_1},x))\|_X&\overset{\eqref{abOSS}}{\leq}\beta(r_2,t)+\gamma\big(\max_{s\in[0,t]}\|y(s,\phi(T_{r_1},x))\|_Y\big)\\
&\overset{\eqref{eq:beta-ri}}{\leq} \frac{r_2}{4}+\frac{r_2}{4}=r_3.
\end{align*}
Hence,
$\|\phi(t,x)\|_X\le r_3$ for $t\in[T_{r_1}+T_{r_2},T].$

Repeating the preceding argument inductively, we get
\begin{equation*}\|\phi(t,x)\|_X\leq r_{i+1},\quad t\in\bigl[T_{r_1}+\ldots+ T_{r_i}, T\bigr],\, i=1,...,k.\end{equation*} 

Taking $i=k$, we obtain
\begin{equation*}\|\phi(t,x)\|_X\leq r_{k+1}<\varepsilon,\quad t\in[T(r,\varepsilon), T].\end{equation*} 

Thus, statement (ii) of Proposition \ref{prop iff} holds. It follows from Proposition \ref{prop iff} that $\Sigma$ is  UGASMO.

\medskip 

($\Leftarrow$) We now prove the reverse implication. Suppose that $\Sigma$ is UGASMO with some $\rho\in\mathcal{K}_\infty$ and $\theta\in \mathcal{KL}$. Fix $x\in X$ and $t\in [0,t_m(x)).$  We divide the proof into two cases.
    
\textbf{Case 1.}   Suppose that
    \begin{align}\label{eq:case}
    \|\phi(s,x)\|_X\ge\rho(\|y(s,x)\|_Y),\quad s\in[0,t].
     \end{align}
Then, by UGASMO, $\|\phi(s,x)\|_X\le \theta(\|x\|_X,s)$ for all $s\in[0,t].$ In particular,
$\|\phi(t,x)\|_X\le \theta(\|x\|_X,t).$

    \textbf{Case 2.}  Now suppose that the condition \eqref{eq:case} does not hold on the whole interval $[0,t].$ Define
    \begin{equation*}A:= \{s\in[0,t]\mid \|\phi(s,x)\|_X< \rho(\|y(s,x)\|_Y)\}\end{equation*}
and let $\tau:=\sup A$.

    If $\tau=t$, continuity of $t\mapsto \phi(t,x)$ and $t\mapsto y(t,x)$ yields
    \begin{align}\label{eq:case2.1}
    \|\phi(t,x)\|_X\le \rho(\|y(t,x)\|_Y)\le \rho(\max_{s\in[0,t]}\|y(s,x)\|_Y).
 \end{align}

If $\tau<t$, then, by the definition of $\tau$ and continuity,
\begin{equation*}\|\phi(s,x)\|_X\ge \rho(\|y(s,x)\|_Y),\quad s\in[\tau,t].\end{equation*}

Moreover, continuity implies
 \begin{equation*}\|\phi(\tau,x)\|_X\le \rho(\|y(\tau,x)\|_Y)\le \rho(\max_{s\in[0,t]}\|y(s,x)\|_Y).\end{equation*}

Furthermore, by the semigroup property and $y(s,z)=h(\phi(s,z))$  for all $s\in[0,t-\tau]$, we have
\begin{align*}
    y(s, \phi(\tau,x))&=h(\phi(\tau+s,x))
    =y(\tau+s,x).
\end{align*}

Consequently, we obtain
\begin{align*}
   \|\phi (s, \phi(\tau,x))\|_X\ge \rho(\|y(s, \phi(\tau,x))\|_Y),\quad s\in[0,t-\tau].
\end{align*}

 UGASMO  applied to the trajectory  from $\phi(\tau,x)$ gives
\begin{align}\label{eq:case2.2}
    \|\phi(t,x)\|_X&=\|\phi(t-\tau,\phi(\tau,x))\|_X\notag\\
    &\le \theta (\|\phi(\tau,x)\|_X,t-\tau)\notag\\
    &\le \theta(\rho(\max_{s\in[0,t]}\|y(s,x)\|_Y),t-\tau)\notag\\
    &\le \theta(\rho(\max_{s\in[0,t]}\|y(s,x)\|_Y),0).
\end{align}

Combining Case 1, \eqref{eq:case2.1} and \eqref{eq:case2.2}, we get
\begin{align*}
\|\phi(t,x)\|_X\le\;&\theta(\|x\|_X,t)+\rho(\max_{s\in[0,t]}\|y(s,x)\|_Y)\\
     &+\theta(\rho(\max_{s\in[0,t]}\|y(s,x)\|_Y),0).
\end{align*}
Therefore, $\Sigma$ is OSS with $\beta(r,t):=\theta(r,t)$ and $\gamma(r):=\rho(r)+\theta(\rho(r),0)$.
\end{proof}

\begin{remark}
    For infinite-dimensional linear systems over Banach spaces with bounded output operators, \cite[Proposition 3.4]{chen2026detectability} shows that exponential OSS with linear gain implies uniform global exponential stability modulo output (UGESMO) via a direct contraction argument on successive time intervals. For general systems, we prove that OSS implies UGASMO using an equivalent  UGASMO characterization.
\end{remark}

\section{Converse OSS Lyapunov theorem}\label{sec:converse OSS}

In this section, we establish a converse OSS Lyapunov theorem. Together with Theorem~\ref{thm:oss-equi-ugasmo}, this yields the following characterization of OSS.

\begin{theorem}[Lyapunov criterion for OSS] 
\label{thm:criterion-OSS}
Let $\Sigma$ satisfy the BIC property.  The following are equivalent:
     \begin{enumerate}[(i)]        
\item  $\Sigma$ admits a coercive OSS Lyapunov function in implication form.

\item  $\Sigma$ is OSS.

\item  $\Sigma$ is UGASMO.
\end{enumerate}
\end{theorem}

\begin{proof}
The implication (i) $\Rightarrow$ (ii) follows from Theorem~\ref{thm:abdirect OSS LT}. The equivalence of (ii) and (iii) is a consequence of Theorem~\ref{thm:oss-equi-ugasmo}.   
The implication (iii) $\Rightarrow$ (i) will be shown later in Theorem~\ref{prop:E0V1} under the BIC property.
\end{proof}

\begin{remark}
    For infinite-dimensional linear systems over Banach spaces with bounded output operators, the corresponding exponential properties are also equivalent: the existence of a coercive exponential OSS Lyapunov function in implication form, exponential OSS, and UGESMO; see \cite[Theorem 3.3]{chen2026detectability}.
\end{remark}

In the remainder of this section, we will prove the implication (iii) $\Rightarrow$ (i) of Theorem~\ref{thm:criterion-OSS}. Throughout this section, assume that $\Sigma$ is UGASMO, and let $\rho$ be as in Proposition \ref{prop iff}.
We can choose   $\hat{\rho}\in \mathcal{K}_{\infty} $  smooth on $(0,\infty)$ with $\hat{\rho}\geq\rho$; see \cite[Proposition A.16]{mironchenko2023}. Thus,  without loss of generality, we may assume that $\rho$ is smooth on $(0,\infty)$.
Define the  sets 
\begin{align*}
\mathcal{D}_{\rho}&:=\{x\in X\mid\|x\|_X\leq\rho(\|h(x)\|_{Y})\},\\
\mathcal{E}_\rho&:=X\setminus\mathcal{D}_{\rho}=\{x\in X\mid \|x\|_X>\rho(\|h(x)\|_{Y})\}.
\end{align*}
Note that $\mathcal{D}_{\rho}$ is always nonempty since $0\in\mathcal{D}_{\rho}$.

\begin{remark}\label{close}
    In view of the continuity of $ \rho$ and $h
    $, $\mathcal{D}_{\rho}$ is a closed set. 
\end{remark}

For every $x\in X$, define the first entrance time into the set $\mathcal{D}_{\rho}$ by
\begin{align}
\label{eq:lambda-map}
\lambda(x):=\inf\{t\in[0,t_m(x))\mid \phi(t,x)\in\mathcal{D}_{\rho}\},
\end{align}
with the convention that $\lambda(x):=t_m(x)$ if  $\phi(t,x)\notin\mathcal{D}_{\rho} $ for all $t\in[0,t_m(x))$,
 which means that the trajectory never enters $\mathcal{D}_{\rho}$.

Let us first establish some basic properties of  $\lambda$.
\begin{lemma}
\label{lem:Value-at-lambda}
For every $x \in X$, we have:
\begin{enumerate}[(i)]
    \item If 
   $\lambda(x)<t_m(x)$, then $\phi(\lambda(x),x)\in \Dc_\rho.$
    \item Assume that $\Sigma$ satisfies the BIC property. If $\lambda(x)=t_m(x),$ then $t_m(x)=\infty.$

\end{enumerate}
 
\end{lemma}

\begin{proof}
 (i).  This is immediate, since  $\mathcal{D}_\rho$ is closed (Remark \ref{close})  and trajectories are continuous.

    (ii). Now suppose that $\lambda(x)=t_m(x).$ Then $\phi(t,x)\in\mathcal{E}_\rho$ for all $t\in[0,t_m(x)).$ The UGASMO estimate \eqref{eq:theta} yields
\begin{equation*}\sup_{t\in[0,t_m(x))}\|\phi(t,x)\|_X\le \theta(\|x\|_X,0)<\infty. \end{equation*}
By the BIC property, it follows that $t_m(x)=\infty.$
\end{proof}

\begin{figure}[!ht]
\centering
\resizebox{\linewidth}{!}{
\begin{tikzpicture}[scale=2.5, every node/.style={font=\small}]

    \coordinate (X) at (0.2,0.3);
    \coordinate (P1) at (1.1,1.1);
    \coordinate (P3) at (3.3,0.39);
 \path[name path=curve] (X) .. controls (P1) .. (P3);
    \path[name path=ellipse] (2.8, 0.1) ellipse (1.2cm and 0.5cm);

    \path[name intersections={of=curve and ellipse, by=I}];
\path[fill=green!8]

  plot[smooth cycle, tension=0.9] coordinates {
    (-0.25,0.85)
    (0.55,1.35)
    (1.75,1.35)
    (3.85,0.95)
    (4.05,0.10)
    (3.70,-0.45)
    (2.30,-0.62)
    (1.45,-0.58)
    (-0.1,-0.48)
    (-0.35,0.05)
  };
  \draw[thick, fill=red!15] (2.8, 0.1) ellipse (1.2cm and 0.5cm);
    \draw[thick,
        postaction={decorate},
        decoration={markings, mark=at position 0.3 with {\arrow{>}}}
    ] (X) .. controls (P1) .. (I);

    \draw[thick] (I) to[bend left=5] (P3);
   \fill (X) circle (0.02);

    \fill (I) circle (0.02);
   \node at (0,1.35) {$X=\mathcal{E}_\rho \cup \mathcal{D}_\rho$};
    \node[below left] at (X) {$x$};
    \node[below=1.2cm,  xshift=1.5cm, overlay] at (X) {$\mathcal{E}_\rho = \{ x \in X \mid \|x\|_X > \rho(\|h(x)\|_Y) \}$};
    \node[above right, xshift=-2.4cm, yshift=-1cm]  at (P1) {$\phi(\cdot, x)$};
    \node[above right] at (I) {$\phi(\lambda(x), x)$};
    \node at (2.8, 0.1) {\( \mathcal{D}_{\rho} = \{ x \in X \mid \|x\|_X \leq \rho(\|h(x)\|_Y) \} \)};

\end{tikzpicture}}
 \caption{Trajectory from $x$ to $\mathcal{D}_{\rho}$ at time $\lambda(x)$.}
\label{entrance time}
\end{figure}
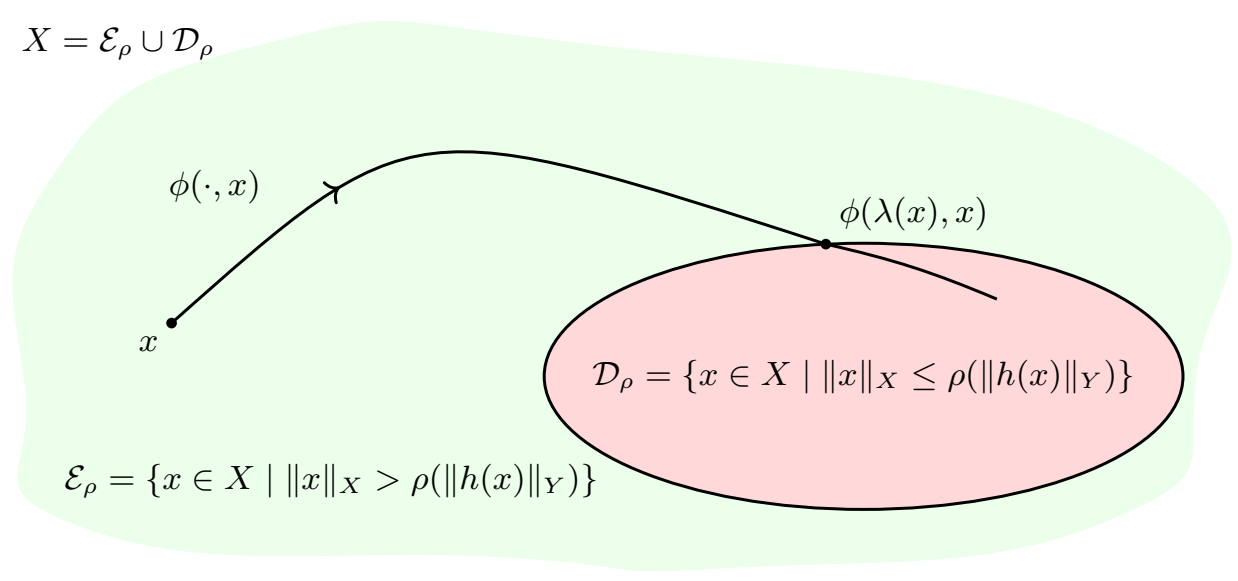

Figure~\ref{entrance time} illustrates the definition of $\lambda$ and Lemma \ref{lem:Value-at-lambda}.

The following lemma describes how the entrance time $\lambda$ changes along trajectories.

\begin{lemma}\label{le:shift-lambda}
Assume that $\Sigma$ satisfies the BIC property.  Then, for $x\in \mathcal{ E}_\rho$ and $t\in[0,\lambda(x))$, we have $
\lambda(\phi(t,x))=\lambda(x)-t
 $, where $\infty-t:=\infty.$
\end{lemma}

\begin{proof}
 Let $x\in \mathcal{ E}_\rho$ and $t\in[0,\lambda(x)).$ If $\lambda(x)<\infty$, then  Lemma~\ref{lem:Value-at-lambda}(ii) and the BIC property imply that $\lambda(x)<t_m(x)$. Thus, by Lemma~\ref{lem:Value-at-lambda}(i), we have $\phi(\lambda(x),x)\in \Dc_\rho$. By the semigroup property, 
$\phi(\lambda(x),x)=\phi(\lambda(x)-t,\phi(t,x))\in \Dc_\rho.$
Moreover, for every $s\in[0,\lambda(x)-t)$,   $\phi(s,\phi(t,x))=\phi(s+t,x)\notin \Dc_\rho.$ 

Consequently,  by the definition of $\lambda$, $\lambda(\phi(t,x))= \lambda(x)-t$ for all $t\in[0,\lambda(x))$.
 If  $\lambda(x)=\infty$, then $t_m(x)=\infty$ by the definition of $\lambda$. Hence, for every $s\ge0$, 
$\phi(s,\phi(t,x))=\phi(s+t,x)\notin \Dc_\rho.$
 Thus  $\lambda(\phi(t,x))=\infty= \lambda(x)-t$.
\end{proof}

Now we can formulate our main result.

\begin{theorem}[Converse OSS Lyapunov theorem]\label{prop:E0V1}
Let $\Sigma$  be a dynamical system satisfying the BIC property. If $\Sigma$ is UGASMO, then it admits a coercive OSS Lyapunov function in implication form.
\end{theorem}

\begin{proof}
Starting from UGASMO, we construct a Lyapunov-like function along the original flow as a weighted supremum taken before the first entrance into the output-dominated region. The semigroup property yields an exponential decay estimate for this function before the trajectory enters $\mathcal{D}_\rho$. A continuous scalar cutoff then extends the function to the whole state space, thereby providing a coercive OSS Lyapunov function in implication form.

Let $\Sigma$ satisfy the BIC property and be UGASMO with $\rho\in \mathcal{K}_\infty$ and $\theta\in \mathcal{KL}$. First, let $x\in\mathcal{E}_\rho$. By the definition of $\lambda(x)$,
\begin{equation*}\phi(t,x)\in\mathcal{E}_\rho,\quad  t\in[0,\lambda(x)).\end{equation*}

Hence, the UGASMO property yields
\begin{equation}\label{eq:phi_theta}
\|\phi(t,x)\|_X\leq \theta(\|x\|_X,t),\quad t\in[0,\lambda(x)).
\end{equation}

Due to Sontag's $\KL$-lemma \cite[Proposition 7]{sontag1998comments}, there exist $\nu_1, \nu_2\in \mathcal{K}_\infty$ such that
\begin{align}\label{eq:theta-estimate}
\theta(r,t)\leq \nu_2(\nu_1(r)e^{-t}),\quad r\geq0, \,   t\geq 0.
\end{align}

Pick  any  function $\omega\in\mathcal{K}_\infty $ such that 
\begin{align}\label{def:omega}
\omega(r)\leq \min\{\nu_1(r),\nu_2^{-1}(r)\},\quad r \ge0.
\end{align}

\vspace{1mm}\noindent\emph{Step 1: Construction of a Lyapunov-like function $V_\rho$, satisfying the dissipation properties on $\mathcal E_\rho$.}

Fix any $\vartheta\in(0,1)$. As $\lambda(x)>0$ for all $x\in \mathcal{E}_{\rho}$, the following map is well-defined:
\begin{align}\label{eq:V1}
V_{\rho}(x):=\sup_{s\in[0,\lambda(x))} e^{\vartheta s} \omega\left(\|\phi(s,x)\|_X\right),\quad x\in \mathcal{E}_{\rho}.
\end{align} 
For all $x\in \mathcal{E}_{\rho}$, the sandwich estimates follow by
\begin{align}\label{eq:V_1coercive}
\omega\left(\|x\|_X\right)\leq V_{\rho}(x)&\overset{\eqref{def:omega}}{\leq }\sup_{s\in[0,\lambda(x))} e^{\vartheta s} \nu_2^{-1}\left(\|\phi(s,x)\|_X\right)\notag\\
&\overset{\eqref{eq:phi_theta}, \eqref{eq:theta-estimate}}{\leq} \sup_{s\in[0,\lambda(x))} e^{\vartheta s} \nu_1\left(\|x\|_X\right)e^{-s}\notag\\
&\leq \nu_1\left(\|x\|_X\right).
\end{align}

For all $x\in \mathcal{E}_{\rho}$ and all $\tau\in(0,\lambda(x))$, by Lemma~\ref{le:shift-lambda},
\begin{align}\label{eq:V-tau}
V_{\rho}(\phi(\tau,x))&=\sup_{s\in[0,\lambda(\phi(\tau,x)))} e^{\vartheta s} \omega\left(\|\phi(s,\phi(\tau,x))\|_X\right)\notag\\
&=\sup_{s\in[0,\lambda(x)-\tau)} e^{\vartheta s} \omega\left(\|\phi(s+\tau,x)\|_X\right)\notag\\
&=\sup_{s\in[\tau,\lambda(x))} e^{\vartheta (s-\tau)} \omega\left(\|\phi(s,x)\|_X\right)\notag\\
&=e^{-\vartheta \tau}\sup_{s\in[\tau,\lambda(x))} e^{\vartheta s} \omega\left(\|\phi(s,x)\|_X\right)\notag\\
&\leq e^{-\vartheta \tau}V_{\rho}(x).
\end{align}
By the latter estimate, for all $x\in \mathcal{ E}_{\rho}$, $t\mapsto V_{\rho}(\phi(t,x))$ is nonincreasing on  $[0,\lambda(x))$. 

In view of the previous computations, to prove the continuity of $t \mapsto V_{\rho}(\phi(t,x))$ on $[0,\lambda(x))$, it is sufficient to prove that, for fixed $x\in\mathcal{E}_\rho$, the function 
\begin{equation*}
t\mapsto \psi_x(t):=\sup_{s\in[t,\lambda(x))} e^{\vartheta s}\omega\left(\|\phi(s,x)\|_X\right)
\end{equation*}
is continuous on $[0,\lambda(x)).$

If $\lambda(x)<\infty$, then Lemma~\ref{lem:Value-at-lambda}(ii), which relies on  the BIC property, implies that $\lambda(x)<t_m(x)$. Hence, by Lemma~\ref{lem:Value-at-lambda}(i), $\psi_x(t)=\max_{s\in[t,\lambda(x)]} e^{\vartheta s}\omega\left(\|\phi(s,x)\|_X\right)$.

If $\lambda(x)=\infty,$  then $t_m(x)=\infty$ by the definition of $\lambda$. By \eqref{eq:theta}, \eqref{eq:theta-estimate} and \eqref{def:omega}, we have for all $s\ge0$
\begin{align*}
0\le e^{\vartheta s}\omega\left(\|\phi(s,x)\|_X\right)
&\overset{\eqref{def:omega}}{\le} e^{\vartheta s} \nu_2^{-1}(\|\phi(s,x)\|_X)\\
&\stackrel{\eqref{eq:theta},\eqref{eq:theta-estimate}}{\le} e^{\vartheta s}  \nu_2^{-1}(\nu_2(\nu_1(\|x\|_X)e^{-s})) \\
&\le \nu_1(\|x\|_X) e^{-(1-\vartheta)s}.
\end{align*}
Thus, the function
$s\mapsto e^{\vartheta s}\omega(\|\phi(s,x)\|_X)$
extends continuously to $[0,\infty]$ by setting its value at infinity to zero. In both cases, the application $t\mapsto[t,\lambda(x)]$ is continuous with compact values. It then follows from Berge's maximum theorem \cite[Theorem~17.31]{AliprantisBorder2006} that $\psi_x$ is continuous.

Moreover, by \eqref{eq:V_1coercive} and \eqref{eq:V-tau}, it holds that
\begin{align}\label{LieV1}
\dot{V}_{\rho}(x)
\overset{\eqref{eq:V-tau}}{\leq }\mathop{\overline{\lim}}\limits_{\tau\rightarrow+0}\frac{e^{-\vartheta \tau}-1}{{\tau}}V_{\rho}(x)
&= -\vartheta V_{\rho}(x)\notag\\
&\overset{\eqref{eq:V_1coercive}}{\leq} -\vartheta \omega\left(\|x\|_X\right).
\end{align}

\vspace{1mm}\noindent\emph{Step 2: Extension of $V_\rho$.}

We now extend $V_{\rho}$ to the whole space $X.$

To construct an OSS Lyapunov function, we introduce a continuous  function depending on the relative size of the state and the output.

Let $\eta: [0,1]\to [0,1]$ be continuous such that
\begin{align*}
    \eta(s)=0,\quad 0\le s\le \frac 23\ ~\text{and}\ ~\eta(s)=1, \quad \frac{3}{4}\le s\le 1.
\end{align*}

Moreover, for $x\in X \setminus \{0\}$, define
\begin{align*}
    q(x):= \frac{\|x\|_X}{\|x\|_X+\rho(\|h(x)\|_Y)}.
\end{align*}

For $x\in X \setminus \{0\}$, set
\[
V(x) \!=\! 
\begin{cases}
\omega(\|x\|_X), &\!\!\!q(x)\le \frac23, \\
\eta( q(x) )V_{\rho}(x)+(1-\eta(q(x)))\omega(\|x\|_X),&\!\!\!q(x)> \frac23,
\end{cases}
\]
and define $V(0):=0.$ The map $V$ is well-defined. Indeed, if $q(x)>\frac23,$ then $\|x\|_X>2\rho(\|h(x)\|_Y)$, and hence $x\in \mathcal{E}_\rho.$ Therefore, $V_\rho(x)$ is defined whenever it appears in the formula for $V.$

We next verify that $V$ is a coercive OSS Lyapunov function in implication form according to Definition~\ref{def:abOSS LF-dissipative}.

\emph{(i): Sandwich bounds.} If $x=0$, then $V(0)=0$. Let $x\in X \setminus \{0\}$. If $q(x)\le \frac23$, then $V(x)=\omega(\|x\|_X)\le \nu_1(\|x\|_X).$ If $q(x)>\frac23,$ then $V(x)$ is a convex combination of $V_\rho(x)$ and $\omega(\|x\|_X),$ and $V_\rho(x)\ge \omega(\|x\|_X).$ Hence,  by \eqref{eq:V_1coercive} and
$\eta(r)\in[0,1]$ for all $r\in[0,1]$,
\begin{align}\label{eq:V-coercivity}
    \omega(\|x\|_X)\le V(x)\le  \omega(\|x\|_X)+\nu_1(\|x\|_X),\quad x\in X.
\end{align}

\emph{(ii): Continuity along trajectories.}  Fix $x\in X$ and let $t^*\in[0,t_m(x))$.

First suppose that $\phi(t^*,x)=0.$ Since $\phi(t,x)\to\phi(t^*,x)$ as $t\to t^*$ and \eqref{eq:V-coercivity} holds, it follows that $V(\phi(t,x))\to 0=V(\phi(t^*,x))$, as $t\to t^*$.

Now suppose that $\phi(t^*,x)\neq0.$ Then $q$ is continuous in a neighborhood of $\phi(t^*,x).$ There are three cases. 

If $q(\phi(t^*,x))<\frac23,$  then  $q(\phi(t,x))<\frac23$ for all $t$ sufficiently close to $t^*.$ Therefore $V(\phi(t,x))=\omega(\|\phi(t,x)\|_X)$ in a neighborhood of $t^*$, and the continuity follows.

If $q(\phi(t^*,x))>\frac23,$  then for all $t$ sufficiently close to $t^*$ we have  $q(\phi(t,x))>\frac23$ and thus $\phi(t,x)\in \mathcal{E}_\rho$.
\begin{align*}
V(\phi(t,x)) =\, &\eta( q(\phi(t,x)) )V_{\rho}(\phi(t,x))\\
&+(1-\eta(q(\phi(t,x))))\omega(\|\phi(t,x)\|_X).
\end{align*}
For $t^*>0$, choose $0\le \hat{t}<t^*$ so that the trajectory remains in $\mathcal{E}_\rho$ from $\hat{t}$ until after $t^*$. The semigroup property gives $\phi(t,x)=\phi(t-\hat{t}, \phi(\hat{t},x))$ with $t-\hat{t}<\lambda(\phi(\hat{t},x))$ in a neighborhood of $t^*$. Applying Step 1 to the restarted trajectory with initial state $\phi(\hat{t},x)$ therefore proves continuity of $V_\rho(\phi(t,x))$ in a neighborhood of $t^*$. At $t^*=0$, the original Step 1 argument applies directly to $x$. The remaining terms are continuous, so $t\mapsto V(\phi(t,x))$ is continuous at $t^*.$

The following argument is adapted from the method used in the proof of \cite[Theorem~3.5]{chen2026detectability}.

Suppose that $q(\phi(t^*,x))=\frac23.$ Then $\eta(q(\phi(t^*,x)))=0$, and hence $V(\phi(t^*,x))=\omega(\|\phi(t^*,x)\|_X)$.  We prove continuity by a sequential argument. Let $(t_n)_{n\in \N}\subset [0,t_m(x))$ be arbitrary with $t_n\to t^*.$ Define
\begin{equation*}I_0:=\{n\in\N \mid q(\phi(t_n,x))\le\frac23 \}\ ~ \text{and}\ ~ I_1:=\N\setminus I_0.\end{equation*}

For every subsequence $(t_{n_j})_{j\in \N}$ of $(t_n)_{n\in\N}$ with $n_j\in I_0$ for every $j\in \N$, by continuity of $\omega$, we have
\begin{equation*} V(\phi(t_{n_j},x))=\omega(\|\phi(t_{n_j},x)\|_X) \to V(\phi(t^*,x)),\quad j\to \infty.
\end{equation*}

Assume that $(t_{n_j})_{j\in \N}$ is a subsequence of $(t_n)_{n\in\N}$ with $n_j\in I_1$ for every $j\in\N$. Then
\begin{align*}
V(\phi(t_{n_j},x)) =\, &\eta( q(\phi(t_{n_j},x)) )V_{\rho}(\phi(t_{n_j},x))\\
&+(1-\eta(q(\phi(t_{n_j},x))))\omega(\|\phi(t_{n_j},x)\|_X).
\end{align*}

Since $q(\phi(t_{n_j},x))\to \frac23$ and $\eta$ is continuous with $\eta(\frac23)=0$, we have $\eta(q(\phi(t_{n_j},x)))\to 0$ as $j\to \infty$. Moreover, $\phi(t_{n_j},x)\to \phi(t^*,x)$ as $j\to \infty$, so  $(\phi(t_{n_j},x))_{j\in \N}$ is bounded. Thus, by \eqref{eq:V_1coercive}, the sequence $(V_\rho(\phi(t_{n_j},x)))_{j\in \N}$ is bounded. Therefore
\begin{equation*} \eta(q(\phi(t_{n_j},x))) V_\rho(\phi(t_{n_j},x))\to 0, \quad j\to \infty. \end{equation*}
The remaining term converges to $\omega(\|\phi(t^*,x)\|_X).$
Thus, 
\begin{equation*}V(\phi(t_{n_j},x))\to V(\phi(t^*,x)),\quad j\to \infty.\end{equation*}

As the subsequences corresponding to $I_0$ and $I_1$ both converge to $V(\phi(t^*,x))$, it follows that
\begin{equation*}V(\phi(t_n,x))\to V(\phi(t^*,x)), \quad n\to\infty.\end{equation*}

Since the sequence $(t_n)_{n\in\N}$ was arbitrary, $t\mapsto V(\phi(t,x))$ is continuous at $t^*$. This proves that $V$ is continuous along trajectories.

\emph{(iii): Decay estimate.} Define $\chi(r):=4\rho(r).$ Let $x\in X$ satisfy $\|x\|_X\ge \chi(\|h(x)\|_Y)$.

If $x=0,$ then $V(\phi(t,0))=V(0)=0$, and hence $\dot{V}(0)=0.$ 

Now let $x\neq 0.$ Then $x\in \mathcal{E}_\rho$ and
\begin{align*}
q(x)\ge \frac{\|x\|_X}{\|x\|_X+\frac{1}{4} \|x\|_X}=\frac45>\frac34>\frac23.
\end{align*}

Since $x\neq0$ and $t\mapsto\phi(t,x)$ is continuous, $\phi(t,x)\neq0$ for all sufficiently small $t\ge0$. Hence, by continuity of $t\mapsto q(\phi(t,x)),$ there exists $\tau\in (0,t_m(x))$ such that $q(\phi(t,x))\ge \frac34$ for all $t\in[0,\tau]$. Hence $\eta(q(\phi(t,x)))=1$ on this interval, and therefore
\begin{equation*}V(\phi(t,x)) =V_\rho(\phi(t,x)),\quad t\in [0,\tau].\end{equation*}
Consequently, by \eqref{LieV1},
\begin{equation*}\dot{V}(x)=\dot{V}_\rho(x)\le -\vartheta \omega (\|x\|_X). \end{equation*}

Thus, combined with the case $x=0$, we obtain
\begin{align*}
\|x\|_{X} \geq\chi(\|h(x)\|_{Y})
\quad \Rightarrow\quad  \dot{V}(x)\leq-\tilde{\alpha}(\|x\|_{X}),
\end{align*}
where $\tilde{\alpha}(r):=\vartheta \omega (r)\in\mathcal{K}_\infty.$
\end{proof}
\begin{remark}\label{re:ugasmo-bic}
   The BIC property is used in Step 1 of the proof of Theorem~\ref{prop:E0V1}, where Lemma~\ref{lem:Value-at-lambda}(ii)  ensures that $\lambda(x)<t_m(x)$ whenever $\lambda(x)<\infty$.
\end{remark}

\section{The vanishing output vanishing state property}\label{sec:vovs}
Motivated by the converging-input converging-state   property considered in \cite{Sontag1989}, we introduce the following property  on the maximal interval of existence. We then investigate the relations between VOVS, UGASMO, and zero-detectability.

\begin{definition}
\label{def:VOVS}
System $\Sigma$ has the \emph{vanishing output vanishing state (VOVS) property} if, for all $x\in X$, 
\[
\lim_{t\to t_m(x)}y(t,x) = 0 \quad\Rightarrow\quad \lim_{t\to t_m(x)}\phi(t,x) = 0.
\]

\end{definition}

\begin{remark}\label{re:vovs-zero-dete}
  The VOVS property implies \emph{zero-detectability} in the sense of \cite{krichman2001input} over the maximal interval of existence. More precisely, for every $x\in X$,
    \[y(t,x)\equiv 0,\quad t\in[0,t_m(x)) \quad\Rightarrow\quad \lim_{t\to t_m(x)}\phi(t,x) = 0.\]
\end{remark}

Let us show that UGASMO systems satisfying the BIC property are VOVS.
\begin{proposition}
\label{prop:UGASMO-implies-weak-convergence}
Suppose that system  $\Sigma$ is UGASMO and satisfies the BIC property. Then it has the VOVS property.
\end{proposition}

\begin{proof}
Let $x\in X$ be such that  $\lim_{t\to t_m(x)} y(t,x)= 0$. 
By assumption, for any $\varepsilon>0$ there is $\tau_1 = \tau_1(\varepsilon)\in[0,t_m(x))$ such that $\|y(t,x)\|_Y\leq \varepsilon$ for all $t\in[\tau_1,t_m(x))$.
Define
\[
\tau_2(\varepsilon):=\inf \{t\in[ \tau_1,t_m(x))\mid\|\phi(t,x)\|_X \leq \rho(\|y(t,x)\|_Y)\},
\]
provided that the set over which the infimum is taken is nonempty. Otherwise, we set $\tau_2(\varepsilon):=t_m(x)$. 

We show that $ \lim_{t\to t_m(x)}\phi(t,x) = 0.$
We first establish that the trajectory $\phi(\cdot,x)$ is bounded on its maximal interval of existence. 

Suppose first that $\tau_2(\varepsilon) = t_m(x)$. Then 
\begin{equation*}\|\phi(t,x)\|_X>\rho(\|y(t,x)\|_Y), \quad t\in[\tau_1,t_m(x)).\end{equation*}

 Applying the UGASMO estimate to the trajectory starting at $\phi(\tau_1,x)$, we obtain, for all $t\in[0,t_m(x)-\tau_1)$, 
\begin{equation}\label{eq:28}
\|\phi(t+\tau_1,x)\|_X\le \theta(\|\phi(\tau_1,x)\|_X,0).
\end{equation}

Thus, the trajectory is bounded on $[\tau_1,t_m(x)),$ and, by continuity, on the compact interval $[0,\tau_1].$ Consequently, 
\begin{align*}
    \sup_{t\in[0,t_m(x))}\|\phi(t,x)\|_X<\infty.
\end{align*}

Now suppose that $\tau_2(\varepsilon)<t_m(x).$ By continuity of trajectories and outputs, we have that
\begin{align}
\label{eq:tau-2-ineq}
\|\phi(\tau_2(\varepsilon),x)\|_X \leq \rho(\|y(\tau_2(\varepsilon),x)\|_Y)\le \rho(\varepsilon).
\end{align}
Now pick any $\hat{t} \in(\tau_2(\varepsilon),t_m(x))$.
Then either
\begin{align}
\label{eq:Alternative1}
\|\phi(\hat{t},x)\|_X \leq \rho(\|y(\hat{t},x)\|_Y)\leq \rho(\varepsilon),
\end{align}
or $\|\phi(\hat{t},x)\|_X > \rho(\|y(\hat{t},x)\|_Y)$. In the latter case, by continuity of trajectories and outputs, there are $t_1 \geq\tau_2(\varepsilon)$ and $t_2 \in (t_1,t_m(x)]$ such that $\hat{t}\in(t_1,t_2)$ and $(t_1,t_2)$ is the largest interval in the sense of set inclusion such that
\begin{align}
\label{eq:aux-Prop44}
\|\phi(t,x)\|_X > \rho(\|y(t,x)\|_Y),\quad t \in (t_1,t_2).    
\end{align}
By continuity, $\|\phi(t,x)\|_X \ge \rho(\|y(t,x)\|_Y)$ for all $ t \in [t_1,t_2)$. Arguing as in the derivation of \eqref{eq:28}, now with initial time $t_1$, we have
\begin{align*}
 \|\phi(t+t_1,x)\|_X \leq\theta(\|\phi(t_1,x)\|_X,t),\quad t\in[0,t_2-t_1).
 \end{align*}
If $t_1 >\tau_2(\varepsilon)$, then the maximality of $(t_1,t_2)$ in \eqref{eq:aux-Prop44} and  continuity imply that $\|\phi(t_1,x)\|_X = \rho(\|y(t_1,x)\|_Y)$. 
Otherwise, if $t_1 =\tau_2(\varepsilon)$, then \eqref{eq:tau-2-ineq} holds, and in any of these cases
\[
\|\phi(t_1,x)\|_X \leq \rho(\|y(t_1,x)\|_Y) \leq \rho(\varepsilon),
\]
and we proceed to
\begin{align*}
 \|\phi(t+t_1,x)\|_X \leq\theta(\rho(\varepsilon),t)\leq\theta(\rho(\varepsilon),0),\quad t\in [0,t_2-t_1).
\end{align*}
Together with \eqref{eq:Alternative1}, we obtain that, for $t\in [ \tau_2(\varepsilon),t_m(x))$,
\begin{align}\label{eq:max-phi}
 \|\phi(t,x)\|_X \leq\max\{\theta(\rho(\varepsilon),0),\rho(\varepsilon)\}.
\end{align}
By continuity, the trajectory is bounded on $[0,\tau_2(\varepsilon)].$ Thus,
\begin{align*}
    \sup_{t\in[0,t_m(x))}\|\phi(t,x)\|_X<\infty.
\end{align*}
In either case, the trajectory is bounded on its maximal interval of existence. The BIC property implies that $t_m(x)=\infty.$ 
 Consequently, $\lim_{t\to \infty} y(t,x)=0.$

It remains to establish convergence of the state. If $\tau_2(\varepsilon)=\infty,$ UGASMO and the semigroup property yield $\lim_{t\to \infty} \phi(t,x)=0$.

If  $\tau_2(\varepsilon)<\infty,$  \eqref{eq:max-phi} gives
\begin{align*}
 \|\phi(t,x)\|_X \leq\max\{\theta(\rho(\varepsilon),0),\rho(\varepsilon)\},\quad t\ge \tau_2(\varepsilon).
\end{align*}

Since $\max\{\theta(\rho(\varepsilon),0),\rho(\varepsilon)\}\to0$ as $\varepsilon\to0,$ and $\varepsilon>0$ is arbitrary, we conclude that
$\lim_{t\to \infty} \phi(t,x)=0$.
\end{proof}

 Theorem~\ref{thm:oss-equi-ugasmo} and Proposition \ref{prop:UGASMO-implies-weak-convergence} yield
\begin{corollary}
\label{cor:OSS-has-vanishing-output-vanishing-state}
If  system $\Sigma$ is OSS and satisfies the BIC property, then it has the VOVS property.
\end{corollary}

\subsection{The VOVS property does not imply UGASMO}

The following example shows that the converse of Proposition~\ref{prop:UGASMO-implies-weak-convergence} does not hold even for linear systems.

Let $X=Y=\ell^2$, where $\ell^2$  is the space of square-summable complex sequences equipped with its standard norm.

Consider the linear system
\begin{align}\label{sys:vovs-not-ugasmo}
    \dot{x}(t)=0, \quad y(t)=Cx(t),
\end{align}
where $C\in\mathcal{L}(\ell^2)$ is defined by
\begin{equation*}Ce_n=\frac1n e_n,\quad n\in\N,\end{equation*}
where $(e_n)_{n\in\N}$ is the standard orthonormal basis of $\ell^2$.

Hence, for every $x\in \ell^2$ and every $t\ge0,$
\begin{equation*}\phi(t,x)=x,\quad y(t,x)=Cx.\end{equation*}
In particular, $t_m(x)=\infty$ for every $x\in X.$

\begin{proposition}\label{prop:vovs-not-ugasmo}
    The following statements hold:
    \begin{enumerate}[(i)]
        \item System \eqref{sys:vovs-not-ugasmo} has the VOVS property.
        \item System \eqref{sys:vovs-not-ugasmo} is not UGASMO.
    \end{enumerate}
\end{proposition}
\begin{proof}
    (i). Since $\phi(t,x)=x$ and $y(t,x)=Cx$ for all $t\ge0$, $y(t,x)\to0$ implies $Cx=0$. 
    By the injectivity of $C$, we obtain $x=0$, and hence $\phi(t,x)\to 0$ as $t\to \infty$.
    Therefore, the system has the VOVS property.

    (ii). Suppose that the system \eqref{sys:vovs-not-ugasmo} is UGASMO with some $\rho\in\mathcal{K}_\infty$ and $\theta\in\mathcal{KL}.$
Since $\rho\in\mathcal{K}_\infty$, there exists $n\in\N$ sufficiently large such that $\rho(\frac1n)\le1$. Moreover, since $\theta\in\mathcal{KL}$, there exists $t^*>0$ such that $\theta(1,t^*)<1.$

Take $x=e_n.$ Then, for every $t\ge0$,
\begin{equation*}\|y(t,e_n)\|_{\ell^2}=\|Ce_n\|_{\ell^2}=\|\frac1n e_n\|_{\ell^2}=\frac1n.\end{equation*}

Thus,  for all $t\in[0,t^*]$,
\begin{align*}
    \|\phi(t,e_n)\|_{\ell^2}=\|e_n\|_{\ell^2}=1\ge \rho(\frac1n)=\rho(\|y(t,e_n)\|_{\ell^2}).
\end{align*}

The UGASMO estimate implies
\begin{equation*}
1=\|\phi(t^*,e_n)\|_{\ell^2}\le \theta (\|e_n\|_{\ell^2},t^*)=\theta (1,t^*)<1,
\end{equation*}
which is a contradiction. Therefore, system \eqref{sys:vovs-not-ugasmo} is not UGASMO. 
\end{proof}

\begin{remark}
    This counterexample relies on the infinite-dimensionality of the state space. Since $C$ is injective, VOVS holds. However, the sequence $(e_n)_{n\in\N}$  satisfies $\|e_n\|_X=1$ and $\|Ce_n\|_X\to 0$, which is possible only in infinite dimensions.
\end{remark}

\begin{remark}
   For infinite-dimensional linear systems over Banach spaces, UGESMO, exponential OSS and OSS are equivalent; see \cite[Theorem 4.5]{chen2026detectability}. As UGESMO implies UGASMO,  system \eqref{sys:vovs-not-ugasmo}, which is not UGASMO, has none of these properties.
\end{remark}

\subsection{Zero-detectability does not imply VOVS}

The following example is taken from \cite[Section 5.1]{chen2026detectability}.

Let $X=L^2(0,1)$ be the space of square integrable Lebesgue measurable complex-valued functions on $(0,1)$ and $Y=\mathbb C$. Consider the system
\begin{equation}\label{eq:example system}
\dot x(t)=Ax(t), \quad y(t)=Cx(t),
\end{equation}
where $A\in\mathcal L(X) $ and  $C\in\mathcal L(X,Y)$
 are defined by
\[
(Af)(s)=is f(s), \ \ s\in (0,1); \quad
Cf=\int_0^1 f(s)\, \mathrm ds.\]
The operator $A $ generates the  $C_0$-group  $T(\cdot)$ given by  
\begin{center}
$(T(t)f)(s)=e^{ist}f(s)$,\quad $s\in (0,1)$,\ $t\in\R$.     
\end{center}
Indeed, $T$ is isometric, i.e., $\|T(t)f\|_X = \|f\|_X$ for all $f\in X$ and all $t\geq 0$.  In particular, $t_m(x)=\infty$ for every $x\in X.$

\begin{proposition}\label{prop:zero-dete-not-Vovs}
 The following hold: 
\begin{enumerate}
    \item [(i)]  System \eqref{eq:example system} is  zero-detectable.
    \item [(ii)]  System \eqref{eq:example system} does not have the VOVS property.
\end{enumerate}
\end{proposition}
\begin{proof}
(i).  It was shown in \cite[Proposition~5.1]{chen2026detectability} that system \eqref{eq:example system} is exponentially zero-detectable. Hence, it is zero-detectable, and thus (i) holds.

(ii). Let $f\in L^2(0,1)$, whence  $f\in L^1(0,1)$. Let $\tilde f$ denote the extension of $f$ by zero to $\mathbb R$. Then $\tilde f\in L^1(\mathbb R)$ and 
\begin{align*}
CT(t)f
&=\int_0^1 e^{ist}f(s)\,\mathrm ds=\int_{\mathbb R} e^{ist}\tilde f(s)\,\mathrm ds=\sqrt{2\pi}\,\hat{\tilde f}(-t),
\end{align*}
where $\hat{\tilde f}$ denotes the Fourier transform of $\tilde f$, and the last equality follows from its definition; see \cite[Section~9.1, formula (4)]{Rudin1987}. By \cite[Theorem~9.6]{Rudin1987},
\begin{equation*}
\hat{\tilde f}\in C_0(\mathbb R),
\end{equation*}
where $C_0(\R)$ denotes the space of all complex-valued continuous functions on $\R$ that vanish at infinity.

Consequently,
\begin{equation*}
CT(t)f\to 0,\quad t\to\infty.
\end{equation*}

On the other hand, since $T$ is isometric, for every $f\neq 0$,
\begin{equation*}
\|T(t)f\|_X=\|f\|_X>0,\quad t\ge0.
\end{equation*}
 Hence, $T(t)f$ does not converge to zero in $X$, and system \eqref{eq:example system} does not have the VOVS property.
\end{proof}

\section{Output-to-state stability of a nonlinear reaction-diffusion system}\label{sec:example}
In this section, we illustrate the developed OSS Lyapunov theory by a nonlinear reaction-diffusion system with a nonlinear output map.

The following example is based on the reaction–diffusion equation considered in \cite{karafyllis2019global}. We adapt its well-posedness framework to the present OSS setting.

Let $a>0$. Consider the nonlinear reaction-diffusion equation with output
\begin{subequations}\label{eq:nonlinear}
    \begin{align}
    x_t(z,t)&=ax_{zz}(z,t)+f(z, x(z,t)),\label{eq:nonf}\\
    y(t)&=\int_0^1 |x (z,t)| \,\mathrm{d} z,\label{eq:output}
\end{align}
\end{subequations}
where  $x_t(z,t)$ and $ x_{zz}(z,t)$ denote  the first and second partial derivative of $x(z,t)$ with respect to $t$ and  $ z$, respectively.
The equation is defined for $z\in(0,1),\ t>0$,  and is endowed with
the homogeneous Dirichlet boundary conditions 
\begin{align}\label{eq:Dirichlet-boundary}
   x(0,t)=x(1,t)=0, 
\end{align}
and the initial condition $x(z,0)=x_0(z).$

Throughout this example, we impose the following assumptions on the reaction term.
\begin{assumption}\label{assm:f}
Assume that $f:[0,1]\times \R\to \R$ satisfies the following conditions:
\begin{enumerate}[(i)]
    \item $f$ is continuous and locally Lipschitz continuous with respect to the second argument, uniformly with respect to $z$ on compact subsets, that is, for every compact set $S\subset \R$, there exists $L(S)>0$ such that
    \begin{align*}
        |f(z,r)-f(z,s)|\le L(S) |r-s|, \quad z\in[0,1],\, r,s\in S.
    \end{align*}
    \item The reaction term satisfies
    \begin{align}\label{eq:f0}
        f(z,0)=0,\quad z\in[0,1].
    \end{align}
   Moreover, there are $q\in\R$, $c_1,c_2\ge0$, and $p\ge  2$  so that
    \begin{align}\label{eq:f1}
     \hspace{-4mm}   s(f(z,s)-qs)\le c_1s^2-c_2|s|^p,\quad z\in[0,1],\,s\in\R.
    \end{align}
\end{enumerate}
\end{assumption}

\begin{remark}
    Assumption \ref{assm:f} covers several representative classes of reaction terms.
\begin{enumerate}[(i)]
    \item Linear reaction terms. Here
$f(z,s)=\lambda_1(z)s$ with
 $\lambda_1\in C([0,1])$. 
    \item Power-type reaction terms. Here
    \begin{equation*}f(z,s)=\lambda_2(z)s-\lambda_3(z)|s|^{r-2}s,\end{equation*}
    where $r\ge 2,\lambda_2,\lambda_3\in C([0,1]),$ and $\lambda_3$ is bounded below by a positive constant. The classical Allen–Cahn nonlinearity $f(s)=s-s^3$ corresponds to $\lambda_2\equiv\lambda_3\equiv1$ and $r=4$; see \cite{allen1979microscopic}. 

    \item Bounded nonpolynomial reaction terms.  Here
    \begin{equation*}f(z,s)= - \lambda_4(z)\frac{s}{1+s^2},\end{equation*}
    where $\lambda_4\in C([0,1])$  is nonnegative.
\end{enumerate}
\end{remark}

The state space we choose as
\begin{equation*}X:=\{x\in H^2(0,1)\mid x(0)=x(1)=0\}\end{equation*}
equipped with the norm $\|\cdot\|_X:=\|\cdot\|_{L^2(0,1)}.$ 
Here $H^2(0,1)$ denotes the Sobolev space of real-valued square-integrable functions on $(0,1)$ whose weak derivatives up to order two are also square-integrable.
Let $Y:=\R$, and define the output map $h:X\to Y$ by
\begin{equation*}h(x):=\int_0^1|x(z)|\,\mathrm{d} z.\end{equation*}

\begin{remark}
The space $X$,  with the norm $\|\cdot\|_X=\|\cdot\|_{L^2(0,1)}$, is not complete and hence  not  Banach. However, completeness is not required in Definition~\ref{dynamical system}. Moreover, Theorem~3.5 in \cite{karafyllis2019global} guarantees that every solution starting in $X$ remains in $X$ throughout its maximal interval of existence. 
\end{remark}

Next, we show that system \eqref{eq:nonlinear} with \eqref{eq:Dirichlet-boundary}  induces a dynamical system in the sense of Definition~\ref{dynamical system}.

\begin{lemma}\label{lem:ex-dy}
    Suppose that Assumption \ref{assm:f} holds. Then, for every $x_0\in X,$ there exists  $t_m(x_0)\in(0,\infty]$, and system \eqref{eq:nonlinear} with \eqref{eq:Dirichlet-boundary} admits a unique maximal solution
    \begin{equation*}
   x\in C([0,1]\times[0,t_m(x_0))).
   \end{equation*}
   The map $t\mapsto x(\cdot,t)$ is in  $C^1((0,t_m(x_0)), L^2(0,1))$
    and
   $x(\cdot,t)\in X$ for all $t\in [0,t_m(x_0)).$ 

Consequently, defining $\phi(t,x_0):=x(\cdot,t),$ system \eqref{eq:nonlinear} with \eqref{eq:Dirichlet-boundary}, together with $h$, defines a dynamical system $\Sigma=(X, Y,\phi,h)$ in the sense of Definition \ref{dynamical system} with equilibrium $0$.
     
\end{lemma}

\begin{proof}
  For $x\in C([0,1])$ and $z\in[0,1]$, set
  \begin{equation*}(\overline{F}(x))(z):=f(z, x(z)).\end{equation*}

The continuity assumption on $f$  implies that $\overline {F}: C([0,1])\to C([0,1])$ is continuous. Moreover, \eqref{eq:f0}  yields  $\overline{F}(0)=0.$ Assumption \ref{assm:f}(i) implies that $\overline{F}$ satisfies the inequality (3.1) in \cite{karafyllis2019global}; see also Remark 2.2(d) therein. Therefore, applying Theorem 3.5 in \cite{karafyllis2019global} yields the existence,  uniqueness and regularity assertions.

Let 
\[
E:=\{x\in C([0,1])\mid x(0)=x(1)=0\}
\]
be equipped with the supremum norm.
By Theorem 3.5 in \cite{karafyllis2019global}
and the variation-of-constants representation,
the maximal PDE solutions coincide with the
maximal mild solutions in $E$ and remain in $X$.

Next, we apply Lemma~\ref{lem:semi-evolution}.
Let $A:D(A)\subset E\to E$ be the Dirichlet realization of
$a\partial_{zz}$ on $E$, defined by
\[
Ax=ax_{zz},\quad
D(A):=\{x\in C^2([0,1])\cap E\mid x_{zz}\in E\}.
\]
By \cite[Theorem VI.4.6]{engel2000one},  $A$ generates  an analytic $C_0$-semigroup on $E$.  Moreover,
$F:=\overline{F}|_E:E\to E$ is locally
Lipschitz continuous by Assumption~\ref{assm:f}.
Moreover, (in line with Remark~\ref{rem:X-norm}) we have continuity of the identity map $I:(X,\|\cdot\|_E) \to (X,\|\cdot\|_X)$ because
\begin{equation*}
\|x\|_X\leq\|x\|_E,\quad x\in X.
\end{equation*}

Since $h$ is Lipschitz continuous on $X$,
Lemma~\ref{lem:semi-evolution} shows that the equation defines a dynamical system on $X$.
Finally, \eqref{eq:f0} and uniqueness imply that \eqref{eq:equilibrium} holds.
\end{proof}

We next establish OSS of the system by means of a coercive OSS Lyapunov function.
\begin{proposition}
 Suppose that Assumption \ref{assm:f} holds. Then 
 \begin{align*}
     V(x):=\frac12 \|x\|_X^2,\quad x\in X
 \end{align*}
 is a coercive OSS Lyapunov function in dissipation form for system \eqref{eq:nonlinear} with \eqref{eq:Dirichlet-boundary}. Hence, system \eqref{eq:nonlinear} with \eqref{eq:Dirichlet-boundary} is OSS. 
\end{proposition}

A related result for semilinear evolution equations on Hilbert spaces is given in \cite[Lemma~11.2.5]{Curtain2020}. Since the state space $X$ is not complete, this result is not directly applicable here.

\begin{proof}
    Obviously,  $V$ is coercive and continuous on $X.$

Let $x(\cdot,t)$ be the maximal solution corresponding to $x_0\in X$. By Lemma~\ref{lem:ex-dy},
$t \mapsto x(\cdot,t)$ belongs to $ C^1\bigl((0,t_m(x_0)), L^2(0,1)\bigr)$
and $x(\cdot,t)\in X$ for all $t\in[0,t_m(x_0))$. Thus, for every $t\in(0,t_m(x_0))$, the time derivative is $L^2$-valued, and the spatial $H^2$ regularity, together with the boundary conditions, justifies the integration by parts below; see \cite[Corollary~8.10]{brezis2011functional}. Here, spatial derivatives are understood in the weak sense. For simplicity, we omit the explicit dependence on $t$ whenever there is no confusion. Hence, by \eqref{eq:nonlinear}, \eqref{eq:Dirichlet-boundary}, and integration by parts, for every $t\in(0,t_m(x_0))$ we obtain
    \begin{align}\label{eq:V-d}
       \frac{\mathrm{d}}{\mathrm{d}t} &V(x(\cdot,t))=\int_0^1 x(z)x_t(z) \,\mathrm{d}z\notag\\
        =\,&a \int_0^1 x(z) x_{zz} (z)\,\mathrm{d}z +\int_0^1 x(z) f(z,x(z)) \,\mathrm{d}z\notag\\
        =\,& a[x(z)x_z(z)]_{z=0}^{z=1}-a \int_0^1 x_z^2(z)\,\mathrm{d}z+ \int_0^1 x(z) f(z,x(z)) \,\mathrm{d}z\notag\\
        =\, & -a \int_0^1 x_z^2(z)\,\mathrm{d}z+ \int_0^1 x(z) f(z,x(z)) \,\mathrm{d}z.
    \end{align}

Let $c_3:=\max\{q+c_1,0\}$. By \eqref{eq:f1} and $p\ge2$, we have
\begin{align}\label{eq:sec}
    \frac{\mathrm{d}}{\mathrm{d}t} V(x(\cdot,t))\le\, & -a \int_0^1 x_z^2(z)\,\mathrm{d}z+(q+c_1)\int_0^1 x^2(z)\,\mathrm{d}z\notag\\
    &-c_2\int_0^1 |x(z)|^p\,\mathrm{d}z\notag\\
    \le\, & -a \int_0^1 x_z^2(z)\,\mathrm{d}z+c_3\int_0^1 x^2(z)\,\mathrm{d}z\notag\\
    &-c_2\Bigl(\int_0^1 x^2(z)\,\mathrm{d}z\Bigr)^{\frac p2}.
\end{align}

Since $x(0)=0$, we get 
\begin{equation*}x(z)=\int_0^z x'(s)\,\mathrm{d}s,\quad z\in[0,1].\end{equation*}

By the Cauchy-Schwarz inequality, we have 
\begin{align*}
|x(z)|
\le \int_0^1 |x_z(s)| \,\mathrm{d}s
&\le \Bigl(\int_0^1 1^2 \,\mathrm{d}s \Bigr)^{\frac12}\Bigl(  \int_0^1 x_z^2(s) \,\mathrm{d}s  \Bigr)^{\frac12}\\
&= \|x_z\|_{L^2(0,1)}
\end{align*}
for all $z\in[0,1]$, and hence $\|x\|_{L^\infty(0,1)}\le  \|x_z\|_{L^2(0,1)}$.

By H\"older’s inequality, we obtain
\begin{align}\label{eq:xz-x}
    \|x\|_{L^2(0,1)}^2&\le \|x\|_{L^\infty(0,1)} \|x\|_{L^1(0,1)}\notag\\
    &\le \|x_z\|_{L^2(0,1)}\|x\|_{L^1(0,1)}.
\end{align}

By  \eqref{eq:xz-x} and Young's inequality, we obtain
\begin{align} \label{eq:c4}
    c_3\int_0^1 x^2(z)\,\mathrm{d}z&\le \|x_z\|_{L^2(0,1)} c_3\|x\|_{L^1(0,1)}\notag\\
    &\le \frac{a}{2} \|x_z\|_{L^2(0,1)}^2 +\frac{c_3^2}{2a}\|x\|_{L^1(0,1)}^2.
\end{align}

By a variation of Wirtinger's inequality \cite[Remark B.1]{smyshlyaev2010adaptive} and $x(1)=0$,  we have
\begin{align}\label{eq:wir}
    -\frac a2\int_0^1x_z^2(z)\, \mathrm dz\leq-\frac{a\pi^2}{8}\int_0^1x^2(z)\, \mathrm dz.
\end{align}

Combining \eqref{eq:sec}, \eqref{eq:c4} and \eqref{eq:wir} yields
\begin{align*}
    \frac{\mathrm{d}}{\mathrm{d}t} V(x(\cdot,t))\le\,& -\frac{a\pi^2}{8}\int_0^1x^2(z)\, \mathrm dz -c_2\Bigl(\int_0^1 x^2(z)\,\mathrm{d}z\Bigr)^{\frac p2}\\
     & +\frac{c_3^2}{2a}\Bigl(\int_0^1 |x(z)| \,\mathrm{d}z\Bigr)^2.
\end{align*}

For all $0<\tau<t<t_m(x_0)$, integration over $[\tau,t]$ yields
\begin{align*}
    V(x(\cdot,t))-&V(x(\cdot,\tau))
    \le \int_\tau^t
    \Big(
    -\frac{a\pi^2}{8}\|x(\cdot,r)\|_X^2
    \\&-c_2\|x(\cdot,r)\|_X^p
    +\frac{c_3^2}{2a}\|x(\cdot,r)\|_{L^1(0,1)}^2
    \Big)\,\mathrm{d}r.
\end{align*}

Letting $\tau\to0$ and using the continuity of $x(\cdot,t)$ in $X$ and the continuity of $V$, we obtain
\begin{align*}
    V(x(\cdot,t))-V(x_0)
    \le \int_0^t
    \Big(
    &-\frac{a\pi^2}{8}\|x(\cdot,r)\|_X^2
    -c_2\|x(\cdot,r)\|_X^p\\
    &+\frac{c_3^2}{2a}\|x(\cdot,r)\|_{L^1(0,1)}^2
    \Big)\,\mathrm{d}r.
\end{align*}

Let
\begin{align*}
    \alpha(s):=\frac{a\pi^2}{8}s^2+c_2s^p,
    \qquad
    \sigma(s):=\frac{c_3^2}{2a}s^2.
\end{align*}
If $c_3=0,$ then $\sigma$ can be chosen as any function in $\mathcal{K}_\infty$.

Since $h(x_0)=\|x_0\|_{L^1(0,1)}$, it follows that
\begin{align*}
    V(x(\cdot,t))\!-\!V(x_0)\!
    \le\! \int_0^t
    \!\bigl(
    -\alpha(\|x(\cdot,r)\|_X)
    +\sigma(h(x(\cdot,r)))
    \bigr)\,\mathrm{d}r.
\end{align*}

Dividing both sides by $t>0$, then letting $t\to0^+$, and using the continuity of $x(\cdot,t)$, $h$, $\alpha$, and $\sigma$, we obtain
\begin{align*}
    \dot V(x_0)
    \le
    -\alpha(\|x_0\|_X)
    +\sigma(h(x_0)).
\end{align*}

Thus, $V$ is a coercive OSS Lyapunov function in dissipation form for system \eqref{eq:nonlinear} with \eqref{eq:Dirichlet-boundary}. Consequently, system \eqref{eq:nonlinear} with \eqref{eq:Dirichlet-boundary} is OSS by Theorem \ref{thm:abdirect OSS LT}.
    \end{proof}

\begin{remark}
The present example is related to the linear parabolic system studied in \cite[Section 7]{chen2026detectability}, for which exponential OSS was established under mixed Robin-Dirichlet boundary conditions and
with the linear output map
\begin{align*}
    h_{\mathrm{lin}}(x)
    = \int_0^1 x(z)\,\mathrm{d}z.
\end{align*}
In contrast, the reaction-diffusion equation \eqref{eq:nonf} has homogeneous Dirichlet boundary conditions and a  linear or nonlinear reaction term. Together with the nonlinear output map \eqref{eq:output}, it defines an OSS system.

\end{remark}

\section{Discussion and outlook}\label{sec:discuss}

\textbf{Comparison with classical ODE results.} 
Our proof of a converse Lyapunov theorem retains the basic idea of analyzing trajectories outside of an output-dominated region, proposed in the ODE setting in the seminal paper \cite{sontag1996}.
However, our construction of an OSS Lyapunov function differs substantially from those in \cite{sontag1996}, \cite{chen2026detectability}. The finite-dimensional argument in \cite{sontag1996} exploits the vector field defining the ODE to introduce an auxiliary slowed-down forward complete system, with the slowing factor chosen using the differentiability of the output map, which cannot be generalized to the infinite-dimensional setting in a straightforward way. The subsequent Lyapunov construction further relies on compactness and a smooth mollifier construction, which do not extend to general normed vector spaces. In contrast, the present proof works directly with the abstract flow, assumes only continuity of the output map, and requires only a continuous cutoff function in the final extension step.

In contrast to \cite{sontag1996}, our converse theorem yields a coercive OSS Lyapunov function \emph{in implication form}. This is the appropriate general conclusion in infinite dimensions, since \cite{chen2026detectability} provides a counterexample showing that an OSS Lyapunov function in dissipation form need not exist even for linear infinite-dimensional systems.
This is in line with other converse Lyapunov results for systems depending on inputs, such as the converse ISS Lyapunov theorem in \cite{MiW17c}, and the converse Lyapunov result for boundedness of reachability sets in \cite{Bachmann26}, which are also shown in implication form. This leads to the thought that Lyapunov functions in implication form are more natural in the infinite-dimensional context, and can be \q{upgraded} to Lyapunov functions in dissipation form only under extra regularity assumptions.

\textbf{Relation to linear OSS theory.} 
A converse Lyapunov theorem for OSS of linear infinite-dimensional systems has been proved in 
\cite{chen2026detectability}. Both papers use the representation of the state space as a union of the sets $\cal D_\rho$ and $\cal E_\rho$, and construct Lyapunov functions separately on both regions gluing them subsequently together. However, the way the Lyapunov functions are constructed and connected together differs in both papers. In particular, thanks to the linearity, a 1-homogeneous OSS Lyapunov function for a linear OSS system is constructed in \cite{chen2026detectability}, which is not to be expected in the general nonlinear case. 
In this work, we have related OSS to the VOVS property, which feels more natural in the nonlinear setting. It is possible to prove the exponential versions of our result in the linear case, but we have not followed this path.

\textbf{Implications for time-delay systems.} 
In \cite[Section 8.10]{Chaillet2023}, the development of the input-output-to-state stability (IOSS) theory for nonlinear time-delay systems was identified as one of the key open problems in the robust stability theory of time-delay systems. 
Under standard well-posedness assumptions, autonomous time-delay systems can be represented as dynamical systems on suitable function spaces; see \cite[Section 7.1]{mironchenko2020input} and \cite{brivadis2024forward}. Hence, equipped with a continuous output map, they fall within the abstract framework considered here. Consequently, the OSS converse Lyapunov characterization established in this paper applies to these time-delay systems.
Extending  our results from OSS to the IOSS setting remains an open problem, but we expect that it is possible to handle this problem by combining the methods developed in this work with the recent advances in converse Lyapunov theory for systems with inputs \cite{Bachmann26}.

\appendix
\section*{Appendix}
\titleformat{\section}
  {\normalfont\normalsize\bfseries\raggedright}
  {\thesection}{0.6em}{}
\section{Chain rule for Dini derivatives}
We  establish an auxiliary lemma that will be used in the proof of Proposition~\ref{prop:p to k}.
This seems to be a folklore result, and we do not claim novelty here, but we were not able to find a reference for this result.

\begin{lemma}[Chain rule for Dini derivatives]
\label{le:chain-rule}
    Let $g:[a,b]\to\R$ be continuous, let $t\in[a,b)$, and let $f:\R\to\R$ be differentiable at $g(t)$. Assume that $f'(g(t))>0$. Then 
    \begin{align*}
        D^+(f\circ g)(t)=f'(g(t)) D^+g(t).
    \end{align*}
\end{lemma}

\begin{proof}
Set $c:=f'(g(t))>0$, and define, for all sufficiently small $\tau>0$,
\[
 q(\tau):=
\begin{cases}
   \frac{f(g(t+\tau))-f(g(t))}{g(t+\tau)-g(t)},\quad &g(t+\tau)\neq g(t),\\
   c, &g(t+\tau)= g(t).
\end{cases}
\]
Then, for all $0<\tau<b-t$,
\begin{align*}
    \frac{f(g(t+\tau))-f(g(t))}{\tau}=q(\tau) \frac{g(t+\tau)-g(t)}{\tau}.
\end{align*}

Since $g(t+\tau)\to g(t)$ as $\tau\to0$ and $f$ is differentiable at $g(t)$, we have $q(\tau)\to c$  as $\tau\to0$.

We use the elementary fact that if $q(\tau)\to c>0$, then (allowing the limits to be equal to $\pm\infty$) 
\begin{align*}
  \mathop{\overline{\lim}}\limits_{\tau\rightarrow+0} q(\tau) k(\tau) =c \mathop{\overline{\lim}}\limits_{\tau\rightarrow+0}  k(\tau)
\end{align*}
for every real-valued function $k$ defined for all sufficiently small $\tau>0.$

Therefore, we obtain
\begin{align*}
    D^+(f\circ g)(t)&=  \mathop{\overline{\lim}}\limits_{\tau\rightarrow+0} q(\tau) \frac{g(t+\tau)-g(t)}{\tau}\\
    &= c  \mathop{\overline{\lim}}\limits_{\tau\rightarrow+0}  \frac{g(t+\tau)-g(t)}{\tau} = f'(g(t)) D^+g(t).
\end{align*}
\end{proof}

\section{Semilinear evolution equations in normed spaces}

The following lemma provides sufficient conditions under which the solution map of a semilinear evolution equation defines a dynamical system on a possibly incomplete normed state space.

Let $(E,\|\cdot\|_E)$ be a Banach space and let $(X,\|\cdot\|_X)$ be a normed linear subspace of $E$.  Let $Y$ be a normed linear space and let $h:X\to Y$ be continuous. Consider the semilinear evolution system
\begin{equation}\label{eq:semi-evo}
    \dot{x}(t)=Ax(t)+F(x(t)),\quad y(t)=h(x(t)),
\end{equation}
where $A:D(A)\subset E\to E$ generates a $C_0$-semigroup on $E$, and $F: E\to E $ is locally Lipschitz continuous on every bounded ball in $E$.

\begin{lemma}
\label{lem:semi-evolution}

Assume that for every $x_0\in X$, the maximal mild solution $\phi(\cdot,x_0)$ of \eqref{eq:semi-evo} in $E$ remains in $X$ throughout its maximal interval of existence and is continuous with respect to the norm of $X$.

Then the solution map $\phi$ defines a dynamical system $\Sigma=(X,Y,\phi,h)$ in the sense of Definition~\ref{dynamical system}.
\end{lemma}

\begin{proof}
    By \cite[Theorem 6.1.4]{pazy2012semigroups}, system \eqref{eq:semi-evo} admits a unique maximal mild solution for every initial state in $E$. Since the equation is autonomous, uniqueness and maximality imply the identity and semigroup properties of the solution map. For every initial state in $X$, the corresponding solution remains in $X$ and is continuous with respect to the norm of $X$ by assumption.  Together with the continuity of $h$, this establishes the assertion.
\end{proof}
\begin{remark}\label{rem:X-norm}
   Let the invariance assumption in Lemma~\ref{lem:semi-evolution} hold. 
   A sufficient condition for the continuity property required in Lemma~\ref{lem:semi-evolution} is the continuity of the identity map $I:(X,\|\cdot\|_E) \to (X,\|\cdot\|_X)$.
\end{remark}


\bibliographystyle{IEEEtran}
\bibliography{References,Mir_LitList_NoMir,MyPublications}


\end{document}